\documentclass[a4paper,12pt,reqno]{amsart}

\usepackage{amsfonts}
\usepackage{amsmath}
\usepackage{amssymb}
\usepackage{mathrsfs}
\usepackage{mathtools}
\usepackage{hyperref}
\usepackage{graphicx}
\usepackage{comment}
\usepackage{nicefrac,xfrac}
\usepackage[usenames]{color}
\usepackage{csquotes}
\usepackage{microtype}

\hypersetup{
colorlinks=true,linkcolor=blue,citecolor=blue,urlcolor=blue,
pdftitle={On Maz'ya's Problem 3.11: Near-Critical Asymptotics and Radial Sharpness},
pdfauthor={Yerkin Shaimerdenov and Nurgissa Yessirkegenov},
pdfkeywords={Hardy inequalities, curl-free vector fields, sharp constants, near-critical asymptotics, radial sharpness, symmetry breaking}}
\numberwithin{equation}{section}

\theoremstyle{plain}
\newtheorem{thm}{Theorem}[section]
\newtheorem{prop}[thm]{Proposition}
\newtheorem{cor}[thm]{Corollary}
\newtheorem{lemma}[thm]{Lemma}
\theoremstyle{definition}
\newtheorem{rem}[thm]{Remark}

\newcommand{\R}{\mathbb{R}}
\newcommand{\Sone}{\mathbb{S}^{1}}
\newcommand{\dd}{\,\mathrm{d}}
\newcommand{\Cp}{\mathscr{C}_{p}}
\newcommand{\Ipg}{\mathcal{I}_{p,\gamma}}
\newcommand{\Jpg}{\mathcal{J}_{p,\gamma}}

\begin{document}

\title[Near-Critical Asymptotics and Radial Sharpness]
{Near-Critical Asymptotics and Radial Sharpness for a Two-Weight Hardy Inequality}

\author[Y. Shaimerdenov]{Yerkin Shaimerdenov}
\address{
  Yerkin Shaimerdenov:
  \endgraf
  SDU University, Kaskelen, Kazakhstan
  \endgraf
  and
  \endgraf
  Institute of Mathematics and Mathematical Modeling, Almaty, Kazakhstan
}

\author[N. Yessirkegenov]{Nurgissa Yessirkegenov}
\address{
  Nurgissa Yessirkegenov:
  \endgraf
  KIMEP University, Almaty, Kazakhstan
}

\thanks{This research is funded by the Science Committee of the Ministry of Science and Higher Education of the Republic of Kazakhstan (Grant No. AP32723635).} 


\begin{abstract}
In this paper, we study the optimal constant $C(p,\gamma)$ in a planar
two-weight Hardy inequality relating the gradient and
Hessian of a scalar potential, posed by Maz'ya in
\cite[Problem~3.11]{Mazya2018}. For every $1\le p<\infty$,
we prove
\[
\lim_{\gamma\to\gamma_c}
2^{p/2}|\gamma-\gamma_c|^p C(p,\gamma)=1,
\qquad \gamma_c=1-\frac2p.
\]
The proof reduces the curl-free constraint to an angular
minimization problem and yields qualitative rigidity of
almost extremizers for $p>1$.
Moreover, we determine the sharp constant in the radial
class for every $1\le p<\infty$ and $\gamma\ne\gamma_c$,
thereby resolving the radial case of Maz'ya's problem.
We also show that the radial inequality remains valid
at $\gamma=\gamma_c$, with the same sharp constant
\[
C_{\rm rad}(p,\gamma)
=\left[1+\left(\gamma-1+\frac2p\right)^2\right]^{-p/2}.
\]
For general, not necessarily radial, potentials, we
identify explicit parameter regions where this constant
remains sharp and establish two sufficient criteria
for symmetry breaking. The exact transition remains
open for $p\ne2$.
In the quadratic case, we recover the known sharp constant
and derive a decomposition of its deficit.
At the critical exponent, we obtain an exact remainder
identity for the logarithmic Hardy inequality, whose
classical constant remains sharp among radial potentials.
The logarithmic inequality fails at $p=1$.
\end{abstract}

\maketitle

\section{Introduction}

\subsection{Maz'ya's problem and the quadratic benchmark}

Let \(1\le p<\infty\), \(\gamma\in\R\), and let
\[
|D^2\psi|^2
=
\psi_{x_1x_1}^2+2\psi_{x_1x_2}^2+\psi_{x_2x_2}^2
\]
denote the Frobenius norm of the Hessian. We consider
\begin{equation}\label{eq:main-ineq}
\int_{\R^2}|\nabla\psi(x)|^p |x|^{p(\gamma-1)}\,\dd x
\le
C(p,\gamma)
\int_{\R^2}|D^2\psi(x)|^p |x|^{p\gamma}\,\dd x.
\end{equation}
Maz'ya posed the problem of finding the optimal constant \(C(p,\gamma)\) for
\[
\gamma\neq 1-\frac2p,
\]
with the additional condition \(\nabla\psi(0)=0\) when
\(\gamma<1-\frac2p\), for \(\psi\in C_0^\infty(\R^2)\) (see
\cite[Problem~3.11]{Mazya2018}). We work throughout with Maz'ya's original
admissible class, which is defined precisely
in Section~\ref{sec:admissible-class}. For \(p=2\), the optimal
constant is
\begin{equation}\label{eq:Mazya-p2}
C(2,\gamma)
=
\begin{cases}
\dfrac{1+(1-\gamma)^2}
{\gamma^2\bigl(3+(1-\gamma)^2\bigr)},
&
-\sqrt3-1\le\gamma\le\sqrt3-1,
\\[1.2em]
\dfrac{1}{1+\gamma^2},
&
\text{otherwise}.
\end{cases}
\end{equation}
Here \(\gamma=0\) is excluded, as in the original problem.

Writing \(u=\nabla\psi\), inequality \eqref{eq:main-ineq} becomes
\begin{equation}\label{eq:curl-free-form}
\int_{\R^2}|u(x)|^p|x|^{p(\gamma-1)}\,\dd x
\le C(p,\gamma)
\int_{\R^2}|\nabla u(x)|^p|x|^{p\gamma}\,\dd x,
\qquad \operatorname{curl}u=0,
\end{equation}
with \(|\nabla u|=|D^2\psi|\). Thus Maz'ya's problem is naturally a sharp weighted
\(L^p\)-Hardy inequality for curl-free fields.
In two dimensions the fixed rotation \(R_\perp(x_1,x_2)=(-x_2,x_1)\)
maps \(u=\nabla\psi\) to the divergence-free field
\(v=R_\perp\nabla\psi\), preserving both \(|v|=|u|\) and
\(|\nabla v|=|\nabla u|\). Thus the planar divergence-free theorem of
Costin--Maz'ya \cite{CostinMazya2008} is exactly the quadratic case of
\eqref{eq:main-ineq}. This planar rotation equivalence is also recorded in
\cite{HamamotoTakahashi2021MA,HamamotoTakahashi2021JFA}.

\subsection{Related works}

Sharp $L^2$ Hardy--Leray and Rellich--Hardy inequalities
for curl-free vector fields have been studied in
\cite{HamamotoTakahashi2021JFA,HamamotoTakahashi2021MA,HamamotoTakahashi2025}.
More recently, Duong and Nguyen established sharp
weighted $L^2$-CKN inequalities and stability estimates
for curl-free fields and second-order derivatives
\cite{DuongNguyen2026}.
In the critical logarithmic two-dimensional $L^2$
case, it is known that the curl-free constraint does
not improve the classical Hardy coefficient
\cite{HamamotoTakahashi2021Critical}.

Related sharp weighted Caffarelli--Kohn--Nirenberg
inequalities for curl-free vector fields and
second-order scalar derivatives were obtained by
Cazacu, Flynn and Lam \cite{CazacuFlynnLam2023},
while exact identities and stability estimates
were developed by Cazacu, Flynn, Lam and Lu
\cite{CazacuFlynnLamLu2024}.
Do, Lam and Lu \cite{DoLamLu2024} investigated
symmetrization and symmetry breaking for weighted
Hardy--Rellich inequalities, whereas Chen and Tang
\cite{ChenTang2024} established second-order $L^p$
Caffarelli--Kohn--Nirenberg identities and sharp
inequalities. The latter two works involve Laplacian
energies, in contrast to the Frobenius norm of the
full Hessian appearing in \eqref{eq:main-ineq}.
To the best of our knowledge, the existing results
do not provide the sharp general-$p$ near-critical
asymptotic for \eqref{eq:main-ineq} or the radial
sharpness and symmetry-breaking criteria established
in the present paper.

\subsection{Main contributions}

Our principal result describes the sharp asymptotic
behavior of the optimal constant $C(p,\gamma)$ over
all admissible potentials as $\gamma$ approaches
the critical value $\gamma_c=1-2/p$. For every $1\le p<\infty$, including $p=1$, we prove
\begin{equation}\label{eq:intro-asymptotic}
\lim_{\gamma\to\gamma_c}
2^{p/2}|\gamma-\gamma_c|^p C(p,\gamma)=1.
\end{equation}
Its proof identifies the first angular (affine) mode as the optimal angular
cell. For $1<p<\infty$, the rigidity theorem describes the original amplitudes
and the mollified derivative fields in the rescaled cylinder variables. In particular, $C(p,\gamma)$ diverges as
$\gamma\to\gamma_c$, whereas the corresponding
radial constant remains bounded.

We next consider Maz'ya's problem restricted to radial
potentials. For every $1\le p<\infty$, we determine
the corresponding sharp constant:
\begin{equation}\label{eq:intro-radial}
C_{\rm rad}(p,\gamma)
=\left[1+\left(\gamma-1+\frac2p\right)^2\right]^{-p/2}.
\end{equation}
We further identify explicit parameter regions for
$p>1$ in which the radial constant remains optimal
among all admissible potentials, without imposing
radial symmetry. In particular, when $p\ge2$, the sufficient conditions are
\begin{equation}\label{eq:intro-positive-region}
\gamma\le -p-\frac2p-\sqrt{2p^2+8}
\quad\text{or}\quad
\gamma\ge -p-\frac2p+\sqrt{2p^2+8}.
\end{equation}
For $1<p<2$ we give a scalar positive-side criterion and explicit sufficient
bounds on both sides. A local angular perturbation and a pure first-harmonic
test prove strict failure of radial sharpness on the union of their sufficient
regions. These conditions are sufficient but do not determine
the exact transition between radial sharpness and
symmetry breaking when $p\ne2$. At $p=2$, the shared method recovers the known
transition $\gamma=-1\pm\sqrt3$ and gives exact remainder decompositions.

\subsection{Strategy and organization}

Our approach is based on a logarithmic change of variables,
which transforms the weighted inequality into a problem on
the cylinder and makes the curl-free constraint explicit.
For the near-critical asymptotic, the key step is a
constrained angular minimization that identifies the
leading angular mode. Convexity remainders are then used
to obtain qualitative rigidity of almost extremizers.
To investigate radial sharpness and symmetry breaking,
we combine exact deficit identities with suitable
angular test functions. In the quadratic case, these
methods recover the known sharp constant and yield
additional deficit decompositions. Finally, at the
critical exponent, a related convexity argument leads
to a logarithmic Hardy inequality with an exact
remainder identity.

The paper is organized as follows. Section~\ref{SEC 2}
introduces the admissible class and notation, together
with the preliminary estimates and identities used
throughout the paper. Section~\ref{SEC 3} presents
the main results, including the near-critical asymptotic,
radial sharpness and symmetry-breaking criteria,
and the quadratic and logarithmic endpoint cases.
The proofs are given in Section~\ref{SEC 4}.

\section{Preliminaries}
\label{SEC 2}

We introduce the admissible class and notation,
and collect the auxiliary estimates and identities
needed in the proofs of the main results.

\subsection{Admissible class and notation}
\label{sec:admissible-class}

For $1\le p<\infty$ and $\gamma\ne1-2/p$, we retain
the admissible class of \cite[Problem~3.11]{Mazya2018}:
\[
\mathcal A_{p,\gamma}
:=\left\{\psi\in C_0^\infty(\R^2;\R):
\nabla\psi(0)=0\ \text{if }\gamma<1-2/p\right\}.
\]
Set
\[
\Ipg[\psi]:=\int_{\R^2}|\nabla\psi|^p|x|^{p(\gamma-1)}\,\dd x,
\qquad
\Jpg[\psi]:=\int_{\R^2}|D^2\psi|^p|x|^{p\gamma}\,\dd x,
\]
and
\begin{equation}\label{eq:a-b}
a:=\gamma-1+\frac2p,
\qquad b:=a-1=\gamma-2+\frac2p.
\end{equation}
The integrals are understood as extended nonnegative integrals. If
\(\Jpg[\psi]=\infty\), the inequality is automatic for any positive constant.
If \(\Jpg[\psi]<\infty\), then \(\Ipg[\psi]<\infty\) as well.
Indeed, for \(a>0\) this follows from boundedness of \(\nabla\psi\) near the origin.
For \(a<0\), let \(k\ge1\) be the order of the first nonzero Taylor term of
\(\nabla\psi\) at the origin. Its derivative is a nonzero homogeneous
polynomial of degree \(k-1\), so finiteness of \(\Jpg\) forces \(a+k>0\).
The same condition makes \(\Ipg\) finite. If all Taylor coefficients vanish,
\(\nabla\psi\) and its derivatives vanish to every order and both integrals
are finite.

Consequently, the best constant for Maz'ya's original class is
\begin{equation}\label{eq:constant-supremum}
C(p,\gamma)
=\sup_{\substack{\psi\in\mathcal A_{p,\gamma}\\0<\Jpg[\psi]<\infty}}
\frac{\Ipg[\psi]}{\Jpg[\psi]}.
\end{equation}
All quotient calculations and deficit identities below are stated for finite
energy. The resulting inequalities hold on the entire original class by
this extended-integral convention. The radial sharpness sequence has an
annular gradient, while its potential may be a nonzero constant near the
origin, as permitted by the original class.

\subsection{Convexity and one-dimensional tools}

We first collect the convexity estimates used to control
the deficits and the one-dimensional inequality needed
for the near-critical lower bound.

\begin{lemma}\label{lem:Cp}
Let \(1<p<\infty\) and \(\xi,\eta\in\R^m\). Define
\[
\Cp(\xi,\eta)
=
|\xi|^p-|\xi-\eta|^p
-p|\xi-\eta|^{p-2}(\xi-\eta)\cdot\eta.
\]
Then \(\Cp(\xi,\eta)\ge0\). Moreover,
\begin{align}
\Cp(\xi,\eta)
&\ge \frac p2|\xi-\eta|^{p-2}|\eta|^2,
&& p\ge2,
\label{eq:Cp-quadratic-explicit}
\\
\Cp(\xi,\eta)
&\ge 2^{2-p}|\eta|^p,
&& p\ge2,
\label{eq:Cp-coercive-large}
\\
\Cp(\xi,\eta)
&\ge
\frac{p(p-1)}2
\frac{|\eta|^2}{\bigl(|\xi-\eta|+|\eta|\bigr)^{2-p}},
&&1<p<2,
\label{eq:Cp-subquadratic-explicit}
\\
\Cp(\xi,\eta)
&\ge
2^{p-3}p(p-1)\bigl(|\xi|+|\xi-\eta|\bigr)^{p-2}|\eta|^2,
&&1<p<2.
\label{eq:Cp-coercive-small}
\end{align}
Expressions at a simultaneous zero are interpreted as zero.
\end{lemma}

\begin{rem}\label{rem:p-less-2}
For \(1<p<2\), if \(\xi_n\) and \(\xi_n-\eta_n\) are bounded in \(L^p\) and
\(\int\Cp(\xi_n,\eta_n)\to0\), then \(\eta_n\to0\) in \(L^p\). Indeed,
\[
\|\eta_n\|_{L^p}^p
\le
\left[
\int
\bigl(|\xi_n|+|\xi_n-\eta_n|\bigr)^{p-2}|\eta_n|^2
\right]^{p/2}
\left[
\int
\bigl(|\xi_n|+|\xi_n-\eta_n|\bigr)^p
\right]^{(2-p)/2}.
\]
\end{rem}

We also record the effective one-dimensional inequality in a form valid down to
\(p=1\).

\begin{lemma}[One-dimensional amplitude inequality]\label{lem:1d-amplitude}
Let \(1\le p<\infty\), \(\sigma\in\{-1,1\}\), and
\(c\in W^{1,p}(\R;\R^2)\). Then
\begin{equation}\label{eq:1d-amplitude}
\int_\R |c'-\sigma c|^p\,\dd s
\ge
\int_\R |c|^p\,\dd s.
\end{equation}
If \(1<p<\infty\), then in fact
\begin{equation}\label{eq:1d-amplitude-remainder}
\int_\R |c'-\sigma c|^p\,\dd s
=
\int_\R |c|^p\,\dd s
+
\int_\R \Cp(c'-\sigma c,c')\,\dd s.
\end{equation}
\end{lemma}

\subsection{Logarithmic cylinder and curl-free constraint}

The logarithmic change of variables removes the weights
and gives a cylinder formulation of the curl-free constraint.

\begin{lemma}[Cylinder formulation]\label{lem:cylinder}
Let \(1\le p<\infty\), and let
\(\psi\in\mathcal A_{p,\gamma}\) satisfy
\(\Jpg[\psi]<\infty\).
Set
\[
a=\gamma-1+\frac2p,
\qquad r=e^t,
\]
and write
\[
e_r=(\cos\theta,\sin\theta),
\qquad
e_\theta=(-\sin\theta,\cos\theta).
\]
Define, in fixed Cartesian components,
\[
Z(t,\theta)
:=r^a\nabla\psi(r\cos\theta,r\sin\theta).
\]
Then
\begin{equation}\label{eq:cylinder-I}
\Ipg[\psi]
=
\int_{\R\times\Sone}|Z|^p\,\dd t\,\dd\theta
\end{equation}
and
\begin{equation}\label{eq:cylinder-J}
\Jpg[\psi]
=
\int_{\R\times\Sone}
\left(|Z_t-aZ|^2+|Z_\theta|^2\right)^{p/2}
\,\dd t\,\dd\theta.
\end{equation}
Moreover,
\begin{equation}\label{eq:cylinder-constraint}
e_r\cdot Z_\theta
=
e_\theta\cdot(Z_t-aZ).
\end{equation}
\end{lemma}

For $a\ne0$, the cylinder formulas give the useful elementary comparison
\begin{equation}\label{eq:comparison-constants}
(1+a^2)^{-p/2}=C_{\rm rad}(p,\gamma)
\le C(p,\gamma)\le |a|^{-p}.
\end{equation}
Indeed, $\mathcal J_{p,\gamma}\ge\int|Z_t-aZ|^p$,
and the supporting-hyperplane inequality at $-aZ$ shows
$\int|Z_t-aZ|^p\ge|a|^p\int|Z|^p$. The linear derivative
term integrates to zero. This also applies for $p=1$ by a subgradient
argument.

For the lower-bound argument it is useful to work in a precise enlarged constrained
class. For \(\varepsilon>0\), \(\sigma\in\{-1,1\}\), and
\(1\le p<\infty\), define \(\mathscr X_{\varepsilon,\sigma}^p\) to be the set of
periodic vector fields \(W:\R\times\Sone\to\R^2\) such that
\[
W\in L^p,\qquad
W_s-\sigma W\in L^p,\qquad
\varepsilon^{-1}W_\theta\in L^p,
\]
and
\begin{equation}\label{eq:scaled-constraint-class}
e_r\cdot\frac{W_\theta}{\varepsilon}
=
e_\theta\cdot(W_s-\sigma W)
\end{equation}
in the sense of distributions. All derivatives in the definition of
\(\mathscr X_{\varepsilon,\sigma}^p\) are distributional. Since
\(W_s,W_\theta\in L^p\), Fubini's theorem gives
\(W(s,\cdot)\in W^{1,p}(\Sone)\) for almost every \(s\), so angular averages
and the periodic Poincar\'e inequality are legitimate on almost every slice. We shall use
\begin{equation}\label{eq:periodic-poincare}
\|F-\overline F\|_{L^p(\Sone)}
\le C_{P,p}\|F_\theta\|_{L^p(\Sone)},
\qquad
\overline F=\frac1{2\pi}\int_{\Sone}F\,\dd\theta,
\end{equation}
for \(1\le p<\infty\). Set
\begin{equation}\label{eq:m-eps}
m_{p,\sigma}(\varepsilon)
:=
\inf_{\substack{W\in\mathscr X_{\varepsilon,\sigma}^p\\ \|W\|_{L^p}=1}}
\int_{\R\times\Sone}
\left(
|W_s-\sigma W|^2+
\left|\frac{W_\theta}{\varepsilon}\right|^2
\right)^{p/2}\,\dd s\,\dd\theta.
\end{equation}
Every rescaled field arising from Lemma~\ref{lem:cylinder} belongs to this class.

\subsection{Angular test functions}

\label{sec:angular-tests}

With the radial benchmark in \eqref{eq:intro-radial}, we construct
nonradial tests whose energy quotient may be smaller.
The first-harmonic test will also be used in the
near-critical asymptotic.

Recall $
a=\gamma-1+2/p\ne0$, $b=a-1$, $r=e^t$, and write
$
V(t,\theta)=e^{bt}\psi(e^t e_r(\theta))$. Since \(\nabla\psi=\psi_r e_r+r^{-1}\psi_\theta e_\theta\), it follows that
\[
\psi_r=e^{-at}(V_t-bV),\qquad
r^{-1}\psi_\theta=e^{-at}V_\theta,
\qquad
Z=(V_t-bV)e_r+V_\theta e_\theta.
\]
Using \(\dd x=e^{2t}\,\dd t\,\dd\theta\),
\(pa=p(\gamma-1)+2\), and \(p(a+1)=p\gamma+2\), the weights cancel:
\[
\Ipg[\psi]=\int_\R\int_0^{2\pi}
\bigl((V_t-bV)^2+V_\theta^2\bigr)^{p/2}\,\dd\theta\,\dd t.
\]
Differentiating once more gives the corresponding Hessian energy
\[
\begin{aligned}
\Jpg[\psi]=\int_\R\int_0^{2\pi}\Bigl[
&\bigl(V_{tt}-(2b+1)V_t+abV\bigr)^2\\
&+2(V_{t\theta}-aV_\theta)^2
+(V_{\theta\theta}+V_t-bV)^2
\Bigr]^{p/2}\,\dd\theta\,\dd t.
\end{aligned}
\]
Here the three expressions are the entries of \(e^{(a+1)t}D^2\psi\)
in the orthonormal polar basis. The off-diagonal entry occurs twice.

Fix \(0\ne\eta\in C_0^\infty((-1,1))\) and a nonzero real-valued smooth
\(2\pi\)-periodic function \(\phi\). Set
\[
V_T(t,\theta)=\eta(t/T)\phi(\theta),\qquad
\psi_T(e^t,\theta)=e^{-bt}\eta(t/T)\phi(\theta).
\]
These potentials are smooth and supported in \(e^{-T}<r<e^T\), so
both energies are finite and \(\nabla\psi_T(0)=0\). Their gradients
are curl-free because they are gradients of smooth scalar functions.

After \(t=Ts\), each \(t\)-derivative of the cutoff contributes a
factor \(T^{-1}\). The terms containing \(\eta'/T\) or
\(\eta''/T^2\) therefore disappear in the following limits:
\[
\begin{aligned}
\frac{\Ipg[\psi_T]}T
&\longrightarrow
\left(\int_\R|\eta|^p\,\dd s\right)
\int_0^{2\pi}\bigl(b^2\phi^2+(\phi')^2\bigr)^{p/2}\,\dd\theta,\\
\frac{\Jpg[\psi_T]}T
&\longrightarrow
\left(\int_\R|\eta|^p\,\dd s\right)
\int_0^{2\pi}
\bigl(a^2b^2\phi^2+2a^2(\phi')^2+(\phi''-b\phi)^2\bigr)^{p/2}\,\dd\theta.
\end{aligned}
\]
Dominated convergence applies: for \(T\ge1\), both rescaled
integrands are bounded by a constant multiple of
\((|\eta|+|\eta'|+|\eta''|)^p\), supported in a fixed compact interval.
Cancelling the common factor \(\int_\R|\eta|^p\) gives
\(\Jpg[\psi_T]/\Ipg[\psi_T]\to\mathcal Q_{p,a}[\phi]\), where
\begin{equation}\label{eq:angular-quotient}
\mathcal Q_{p,a}[\phi]
:=\frac{\displaystyle\int_0^{2\pi}
\bigl(a^2b^2\phi^2+2a^2(\phi')^2+(\phi''-b\phi)^2\bigr)^{p/2}\,\dd\theta}
{\displaystyle\int_0^{2\pi}
\bigl(b^2\phi^2+(\phi')^2\bigr)^{p/2}\,\dd\theta},
\end{equation}
provided the denominator is positive. Since
\(C(p,\gamma)^{-1}=\inf_\psi\Jpg[\psi]/\Ipg[\psi]\), with the
infimum over admissible nonzero functions of finite energy, we obtain
\begin{equation}\label{eq:angular-test-bound}
C(p,\gamma)^{-1}\le\mathcal Q_{p,a}[\phi].
\end{equation}
In particular, \(\mathcal Q_{p,a}[\phi]<(1+a^2)^{p/2}\) implies
\(C(p,\gamma)>C_{\rm rad}(p,\gamma)\), which we call symmetry breaking.

\begin{prop}[First-harmonic test]\label{prop:first-harmonic}
Let \(1\le p<\infty\), \(a=\gamma-1+2/p\ne0\), and \(b=a-1\). Set
\begin{equation}\label{eq:Ip}
I_p(b):=\int_0^{2\pi}(b^2\cos^2\theta+\sin^2\theta)^{p/2}\,\dd\theta,
\end{equation}
\begin{equation}\label{eq:Jp}
J_p(b):=\int_0^{2\pi}(b^2\cos^2\theta+1+\sin^2\theta)^{p/2}\,\dd\theta.
\end{equation}
Then
\begin{equation}\label{eq:first-harmonic-bound}
C(p,\gamma)\ge\frac{I_p(b)}{|a|^pJ_p(b)}.
\end{equation}
Consequently, strict symmetry breaking holds if
\begin{equation}\label{eq:first-harmonic-breaking}
|a|^pJ_p(b)<(1+a^2)^{p/2}I_p(b).
\end{equation}
\end{prop}

\subsection{Nonlinear cylinder identities}

To prove radial sharpness, we also need lower bounds
for the energy of arbitrary admissible potentials.
We therefore derive an exact identity for the deficit
relative to the radial constant.

Write the field from Lemma~\ref{lem:cylinder} as
\[
Z=P e_r+Qe_\theta,\qquad U=(P,Q),\qquad
P=V_t-bV,\quad Q=V_\theta.
\]
Let \(J(P,Q)=(-Q,P)\), the rotation through \(\pi/2\).
Differentiating the polar basis gives
\[
Z_t-aZ=(P_t-aP)e_r+(Q_t-aQ)e_\theta,\qquad
Z_\theta=(P_\theta-Q)e_r+(Q_\theta+P)e_\theta.
\]
Consequently,
\begin{equation}\label{eq:U-cylinder-J}
\Jpg[\psi]
=\int_{\R\times\Sone}
\left(|U_t-aU|^2+|U_\theta+JU|^2\right)^{p/2}
\,\dd t\,\dd\theta,
\end{equation}
while \(\Ipg[\psi]=\int_{\R\times\Sone}|U|^p\).
Since \(P,Q\) come from the same potential \(V\),
\begin{equation}\label{eq:U-constraint}
P_\theta=V_{t\theta}-bV_\theta=Q_t-bQ.
\end{equation}
This is the curl-free constraint in these variables. Below, all
unqualified integrals are over \(\R\times\Sone\), with measure
\(\dd t\,\dd\theta\).

\begin{lemma}[Nonlinear curl-free identity]\label{lem:KL-identity}
Let \(1<p<\infty\), and let \(U=(P,Q)\) arise from an admissible
potential with \(\Jpg[\psi]<\infty\). Set
\[
A=\int |U|^{p-2}Q^2,\qquad
K=\int |U|^{p-2}(PQ_\theta-QP_\theta),
\]
and
\[
L=\int |U|^{p-4}Q
\left[Q\det(U,U_\theta)+P\det(U,U_t)\right],
\]
where \(\det(U,U_\theta)=PQ_\theta-QP_\theta\) and
\(\det(U,U_t)=PQ_t-QP_t\). Then
\begin{equation}\label{eq:KL-identity}
K+(p-2)L=pbA.
\end{equation}
The integrands are defined to be zero at \(U=0\).
\end{lemma}

\begin{lemma}[Exact deficit identity]\label{lem:exact-deficit}
Under the assumptions above, put
\[
\lambda=\sqrt{1+a^2},\qquad
X=(-aU,JU),\qquad Y=(U_t,U_\theta).
\]
Then
\begin{equation}\label{eq:outer-Cp-deficit}
\Jpg[\psi]-\lambda^p\Ipg[\psi]
=p\lambda^{p-2}K+R_p,\qquad
R_p=\int\Cp(X+Y,Y)\ge0.
\end{equation}
\end{lemma}

\subsection{The angular cell and constrained lower limit}

We now use the cylinder formulation to study the limit
$\gamma\to1-2/p$. The main step is a constrained lower
bound, whose sharp coefficient is determined by the
following angular minimization problem.

\begin{lemma}[Angular cell]\label{lem:angular-cell}
Let \(1\le p<\infty\), let \(q\in\R^2\) be fixed, and suppose
\(Y\in L^p(\Sone;\R^2)\) satisfies
\[
\int_{\Sone}Y\,\dd\theta=0,
\qquad
e_r\cdot Y=e_\theta\cdot q
\quad\text{almost everywhere}.
\]
Thus \(Y\) has zero average, and its radial component is prescribed.
Define
\begin{equation}\label{eq:Ycell}
\mathcal Y[q]
=(e_\theta\cdot q)e_r+(e_r\cdot q)e_\theta.
\end{equation}
Then
\begin{equation}\label{eq:angular-cell-lower}
\int_{\Sone}
\left(|q|^2+|Y|^2\right)^{p/2}\,\dd\theta
\ge 2\pi\,2^{p/2}|q|^p,
\end{equation}
with equality for \(Y=\mathcal Y[q]\).

For \(1<p<\infty\), the exact identity
\begin{equation}\label{eq:angular-cell-remainder}
\begin{aligned}
&\int_{\Sone}
\left(|q|^2+|Y|^2\right)^{p/2}\,\dd\theta
\\
&\quad=
2\pi\,2^{p/2}|q|^p
+\int_{\Sone}
\Cp\bigl((q,Y),(0,Y-\mathcal Y[q])\bigr)\,\dd\theta
\end{aligned}
\end{equation}
holds.
\end{lemma}
The preceding lemma identifies the optimal angular
energy. We next use this bound to obtain the
near-critical lower estimate for the constrained
cylinder problem.

\begin{prop}[Near-critical lower bound]
\label{prop:constrained-liminf}
Let \(1\le p<\infty\) and \(\sigma\in\{-1,1\}\). Then
\begin{equation}\label{eq:constrained-liminf}
\liminf_{\varepsilon\downarrow0}m_{p,\sigma}(\varepsilon)
\ge 2^{p/2}.
\end{equation}
\end{prop}

\subsection{A nonlinear estimate for radial sharpness}

The following estimate will be used to obtain explicit sufficient regions of radial sharpness.

\begin{lemma}[A strengthened estimate for the nonlinear term]
\label{lem:L-control}
Under the assumptions of Lemma~\ref{lem:KL-identity}, let
\(p\ge2\) and set
\[
B:=\int_{\R\times\Sone}
|U|^{p-2}\left(|U_t|^2+|U_\theta|^2\right).
\]
With \(A\) and \(L\) defined there, one has
\begin{equation}\label{eq:strong-L-control}
L^2+(L-bA)^2\le AB.
\end{equation}
\end{lemma}

\section{Main results}
\label{SEC 3}

We state the sharp asymptotic, the rigidity conclusion, the radial and nonradial criteria, and the quadratic and logarithmic endpoint results. Their proofs are given in Section~\ref{sec:all-proofs}.

\subsection{The sharp near-critical asymptotic}

Combining Proposition~\ref{prop:constrained-liminf}
with the first-harmonic test of
Proposition~\ref{prop:first-harmonic} gives the
sharp near-critical asymptotic.

\begin{thm}[Near-critical asymptotics]\label{thm:asymptotic}
Let \(1\le p<\infty\). Then
\begin{equation}\label{eq:main-asymptotic}
{
\lim_{\gamma\to1-\frac2p}
2^{p/2}
\left|\gamma-1+\frac2p\right|^p
C(p,\gamma)=1.
}
\end{equation}
Equivalently,
\[
C(p,\gamma)
=
\frac{1+o(1)}
{2^{p/2}\left|\gamma-1+\frac2p\right|^p}
\qquad
\text{as }\gamma\to1-\frac2p.
\]
In particular,
\[
C(1,\gamma)\sim\frac1{\sqrt2\,|\gamma+1|}
\qquad(\gamma\to-1).
\]
\end{thm}

\begin{rem}[Origin of the factor \(2^{p/2}\)]
The factor \(2^{p/2}\) is the energy of the optimal angular cell:
\[
\bigl(|q|^2+|\mathcal Y[q]|^2\bigr)^{p/2}
=2^{p/2}|q|^p.
\]
In the original polar variables this cell corresponds to the first angular, or affine,
mode.
\end{rem}

\subsection{Qualitative rigidity of almost extremizers}

Retaining the convexity remainders in the proof of Proposition~\ref{prop:constrained-liminf} yields a qualitative description
of almost extremizers for $1<p<\infty$.

\begin{thm}[Near-critical rigidity]\label{thm:rigidity}
Let \(1<p<\infty\), and let
\[
a_n=\gamma_n-1+\frac2p\longrightarrow0,
\qquad
\varepsilon_n=|a_n|>0,
\qquad
\sigma_n=\operatorname{sgn}(a_n).
\]
Suppose that \(\psi_n\in\mathcal A_{p,\gamma_n}\) satisfies
\[
\mathcal I_{p,\gamma_n}[\psi_n]=1,
\qquad
\frac{\mathcal J_{p,\gamma_n}[\psi_n]}{\varepsilon_n^p}
\longrightarrow2^{p/2}.
\]

Let \(Z_n\) be the field from Lemma~\ref{lem:cylinder}, and set
\[
W_n(s,\theta)
=\varepsilon_n^{-1/p}Z_n(s/\varepsilon_n,\theta),
\qquad
d_n=(W_n)_s-\sigma_nW_n,
\qquad
Y_n=\varepsilon_n^{-1}(W_n)_\theta.
\]
Fix \(0\le\rho\in C_0^\infty((-1,1))\) with
\(\int_\R\rho\,\dd s=1\), and write
\[
\rho_\delta(s)=\delta^{-1}\rho(s/\delta).
\]
Mollify in \(s\) at scale \(\delta_n=\sqrt{\varepsilon_n}\):
\[
\widetilde W_n=\rho_{\delta_n}*_sW_n,
\qquad
\widetilde d_n=\rho_{\delta_n}*_sd_n,
\qquad
\widetilde Y_n=\rho_{\delta_n}*_sY_n.
\]
Define
\[
c_n(s)=\frac1{2\pi}\int_{\Sone}
\widetilde W_n(s,\theta)\,\dd\theta,
\qquad
d_n^0=c_n'-\sigma_nc_n.
\]

Then the following conclusions hold.

\smallskip
\noindent
\emph{Angular shape and mass.}
The original field approaches its effective amplitude:
\begin{align}
2\pi\int_\R|c_n|^p\,\dd s&\longrightarrow1,
\label{eq:rig-mass}\\
\|W_n-c_n\|_{L^p(\R\times\Sone)}&\longrightarrow0.
\label{eq:rig-W}
\end{align}

\noindent
\emph{Slow variation.}
The amplitude satisfies
\begin{equation}\label{eq:rig-flat}
\|c_n'\|_{L^p(\R)}\longrightarrow0.
\end{equation}

\noindent
\emph{Angular-cell structure.}
The mollified derivative fields satisfy
\begin{align}
\|\widetilde d_n-d_n^0\|_{L^p(\R\times\Sone)}
&\longrightarrow0,
\label{eq:rig-d}\\
\|\widetilde Y_n-\mathcal Y[d_n^0]\|_{L^p(\R\times\Sone)}
&\longrightarrow0.
\label{eq:rig-Y}
\end{align}
In particular,
\[
\|\widetilde d_n+\sigma_nc_n\|_{L^p}\longrightarrow0,
\qquad
\|\widetilde Y_n+\sigma_n\mathcal Y[c_n]\|_{L^p}
\longrightarrow0.
\]
Here and below, unqualified \(L^p\)-norms are over
\(\R\times\Sone\). The derivative conclusions concern the
mollified fields, whereas \eqref{eq:rig-W} concerns \(W_n\)
itself.
\end{thm}

\begin{rem}[Geometric interpretation]
Equation \eqref{eq:rig-W} says that in Cartesian cylinder coordinates an almost
extremizer becomes asymptotically independent of \(\theta\). Taking its polar
components \(U=(e_r\cdot Z,e_\theta\cdot Z)\), this is precisely the first
angular mode in polar coordinates. Thus
the first angular/affine mode is the leading near-critical concentration mechanism.
\end{rem}

\begin{rem}[Slow variation, not a unique plateau]
The conclusions
\[
2\pi\|c_n\|_{L^p}^p\to1,
\qquad
\|c_n'\|_{L^p}\to0
\]
show that the effective amplitude becomes slowly varying and delocalized in the
rescaled logarithmic variable. They do not, by themselves, imply convergence to a
single plateau or exclude several widely separated slowly varying components.
\end{rem}

\begin{rem}[Two rigidity defects]
Set
\[
X_n=(\widetilde d_n,\widetilde Y_n),
\qquad
X_n^0=(d_n^0,\mathcal Y[d_n^0]).
\]
The proof isolates two nonnegative defects,
\[
R_n^{\rm ang}
:=\int_{\R\times\Sone}
\Cp(X_n,X_n-X_n^0)\,\dd s\,\dd\theta
\]
and
\[
R_n^{\rm 1D}
:=\int_\R\Cp(d_n^0,c_n')\,\dd s.
\]

The first measures deviation from the optimal angular cell, while the second measures
deviation from the virtual one-dimensional extremizing mechanism. A genuine
quantitative stability inequality would require turning these sequential conclusions into
an explicit deficit--distance estimate.
\end{rem}

\subsection{The sharp radial constant}

For radial potentials, the cylinder formulation reduces
the problem to one dimension. We compute the sharp
constant, which will serve as the reference value
for the full problem.

\begin{thm}[Sharp radial constant]\label{thm:radial}
Let \(1\le p<\infty\) and \(\gamma\in\R\). Among radial functions
\(\psi(x)=f(|x|)\in C_0^\infty(\R^2)\) for which both weighted
integrals are finite (with the origin condition of the full problem
when applicable), the optimal constant is
\begin{equation}\label{eq:radial-constant}
C_{\rm rad}(p,\gamma)
=
\left[
1+\left(\gamma-1+\frac2p\right)^2
\right]^{-p/2}.
\end{equation}
The radial inequality remains meaningful at \(\gamma=1-2/p\), where
\(C_{\rm rad}=1\).
\end{thm}

\begin{rem}\label{rem:radial-p2}
For \(p=2\), one has \(a=\gamma\), and
\[
C_{\rm rad}(2,\gamma)=\frac1{1+\gamma^2},
\]
which is exactly the second branch in Maz'ya's formula \eqref{eq:Mazya-p2}.
For \(p=1\), Theorem~\ref{thm:radial} gives the explicit endpoint formula
\[
C_{\rm rad}(1,\gamma)=\frac{1}{\sqrt{1+(\gamma+1)^2}}.
\]
\end{rem}

\subsection{Sufficient regions of radial sharpness}

\label{sec:radial-sharpness}

Theorem~\ref{thm:asymptotic} excludes radial sharpness
near the critical exponent. We now return to the
nonlinear deficit identity and the estimate of Lemma~\ref{lem:L-control}
to find regions where the radial constant is sharp.

\begin{thm}[An explicit parameter region of radial sharpness]
\label{thm:radial-sharpness}
Let \(1<p<\infty\) and consider \eqref{eq:main-ineq} for real-valued admissible
functions. Set
\[
a:=\gamma-1+\frac2p,
\qquad
b:=a-1=\gamma-2+\frac2p.
\]
Then
\begin{equation}\label{eq:full-radial-sharp}
{
C(p,\gamma)
=
\left(1+a^2\right)^{-p/2}
}
\end{equation}
under either of the following sufficient conditions.

\begin{enumerate}
\item If \(p\ge2\), it is sufficient that
\begin{equation}\label{eq:radial-region-pge2-b}
{b\le -(p+2)-\sqrt{2p^2+8}\quad\text{or}\quad
b\ge -(p+2)+\sqrt{2p^2+8},}
\end{equation}
or equivalently
\begin{equation}\label{eq:radial-region-pge2-gamma}
{
\gamma\le -p-\frac2p-\sqrt{2p^2+8}\quad\text{or}\quad
\gamma\ge -p-\frac2p+\sqrt{2p^2+8}.
}
\end{equation}

\item If \(1<p<2\), put
\[
\lambda:=\sqrt{1+a^2},
\qquad
\kappa_p:=\frac{p(p-1)}2.
\]
It is sufficient that
\begin{equation}\label{eq:radial-region-sub2}
{
p^2b
\ge
\sup_{z\ge0}
\left\{
p(2-p)\lambda z
-
\kappa_p\lambda^2\frac{z^2}{(1+z)^{2-p}}
\right\}.
}
\end{equation}
\end{enumerate}
\end{thm}

\begin{rem}[The positive-side scalar threshold]\label{rem:implicit-threshold}
For \(1<p<2\), set \(q=2-p\). With \(y=\lambda z\), the supremum in
\eqref{eq:radial-region-sub2} becomes
\[
\sup_{y\ge0}\left\{pqy-\kappa_p\frac{y^2}{(1+y/\lambda)^q}\right\}.
\]
It is finite, continuous in \(\lambda>0\), and nonincreasing in \(\lambda\).
For \(b\ge0\), \(\lambda=\sqrt{1+(b+1)^2}\) increases with \(b\), while
\(p^2b\) strictly increases. Thus the sufficient set in
\eqref{eq:radial-region-sub2} is a ray \([b_*(p),\infty)\).
The next corollary gives an elementary upper bound for \(b_*(p)\), which
can be much larger than the threshold obtained from the scalar supremum.
\end{rem}

\begin{cor}[A simple explicit condition for \(1<p<2\)]
\label{cor:radial-sharpness-sub2}
Let \(1<p<2\) and define
\[
d_p:=2^{p-3}p(p-1),
\]
\[
B_p
:=
\max\left\{
\frac{(2-p)^2}{4d_p},
\frac{p-1}{p^2}(2-p)^{\frac{p}{p-1}}d_p^{-\frac1{p-1}}
\right\}.
\]
If
\begin{equation}\label{eq:radial-region-sub2-explicit}
{
\gamma\ge2-\frac2p+B_p,
}
\end{equation}
then \eqref{eq:full-radial-sharp} holds.
\end{cor}
The preceding result gives an explicit sufficient
condition on the positive side. We now establish
a corresponding condition on the negative side
for $1<p<2$.

\begin{thm}[A negative radial region for \(1<p<2\)]
\label{thm:negative-sub2}
Let \(1<p<2\) and put
\[
D_p=\frac{2^{3-p}\bigl(2p+\sqrt2(2-p)\bigr)}{p-1}.
\]
If
\begin{equation}\label{eq:negative-region-sub2}
\gamma\le 2-\frac2p-D_p,
\end{equation}
then the full constant equals the radial constant
\(C(p,\gamma)=(1+a^2)^{-p/2}\).
\end{thm}

\begin{rem}[Scope of the subquadratic bounds]
Both explicit thresholds deteriorate as \(p\downarrow1\):
\(B_p\to\infty\), and
\[
D_p\sim\frac{4(2+\sqrt2)}{p-1}.
\]
The explicit bound $B_p$ is particularly coarse.
For example, when $p=6/5$, direct evaluation gives
$B_p\approx2.34\times10^4$, whereas numerical
minimization of the scalar criterion
\eqref{eq:radial-region-sub2} yields
$b_*(p)\approx4.034$.
These values illustrate the difference between
the two sufficient bounds and do not determine
the optimal transition between radial sharpness
and symmetry breaking. No outer radial-sharpness region
at \(p=1\) is established here.
\end{rem}

\begin{rem}\label{rem:radial-threshold-not-optimal}
The sufficient thresholds in Theorems~\ref{thm:radial-sharpness} and
\ref{thm:negative-sub2} do not describe the optimal transitions. For example,
when \(p=2\) the general estimate yields \(\gamma\le-7\) or
\(\gamma\ge1\), whereas Corollary~\ref{cor:quadratic-transition} gives the
exact endpoints \(-1\pm\sqrt3\).
\end{rem}

\subsection{Two criteria for symmetry breaking}

\label{sec:symmetry-breaking}

We next identify parameters for which the radial
constant is not sharp. A perturbation of the constant
angular profile gives a second criterion, complementing
Proposition~\ref{prop:first-harmonic}.

\begin{prop}[Local angular instability]\label{prop:instability}
Let \(1<p<\infty\), \(a=\gamma-1+2/p\ne0\), and define
\[
S_p(a):=(a-1)^2+2p(a-1)+\frac{a^2+p-1}{1+a^2}.
\]
If \(S_p(a)<0\), then
\begin{equation}\label{eq:strict-breaking}
C(p,\gamma)>(1+a^2)^{-p/2}.
\end{equation}
\end{prop}

Combining Proposition~\ref{prop:instability} with
\eqref{eq:first-harmonic-breaking} proves strict symmetry breaking on
the union of the two sufficient sets. At \(p=2\), they coincide, since
\(S_2(a)=b^2+4b+1\).

\begin{rem}[Symmetry breaking and the unresolved transition]
The condition \(S_p(a)\ge0\) is necessary for radial sharpness,
not sufficient. For example, at \(p=4\), \(\gamma=5/4\), one has
\(a=3/4\), \(b=-1/4\) and \(S_4(3/4)=137/400>0\), whereas the
first-harmonic test gives
\[
\mathcal Q_{4,3/4}[\cos\theta]
=\frac{37017}{18688}<\frac{625}{256}
=(1+a^2)^2.
\]
Thus the radial constant is not globally sharp there. This comparison
does not establish positivity of the full second variation. In
particular, the exact nonlinear transition cannot be inferred solely
from the sign of \(S_p\).
For $p=4$, the integrals in
\eqref{eq:Ip} and \eqref{eq:Jp} can be
evaluated explicitly:
\[
I_4(b)=\frac{\pi}{4}(3b^4+2b^2+3),
\qquad
J_4(b)=\frac{\pi}{4}(3b^4+10b^2+19).
\]
Consequently, the first-harmonic criterion
\eqref{eq:first-harmonic-breaking} reduces
to the polynomial inequality
\[
(b+1)^4(3b^4+10b^2+19)
<\bigl(1+(b+1)^2\bigr)^2(3b^4+2b^2+3).
\]
These two criteria provide different sufficient
regions for symmetry breaking. Their union therefore
gives a larger sufficient region than either
criterion alone. For $p\ne 2$, the exact value of $C(p,\gamma)$
remains unknown in general. In particular, the precise
transition between radial sharpness and symmetry
breaking remains an open problem.
Theorems~\ref{thm:radial-sharpness} and
\ref{thm:negative-sub2} give sufficient regions
of radial sharpness, while
Propositions~\ref{prop:instability} and
\ref{prop:first-harmonic} provide sufficient
conditions for strict symmetry breaking.
\end{rem}

\subsection{The quadratic case}

For $p=2$, the cylinder identities determine the exact
transition between radial sharpness and symmetry breaking.
We recover the constant of Costin--Maz'ya
\cite{CostinMazya2008} and decompose the deficit into
radial, first angular, and higher angular contributions.
A sharp inequality with a remainder and non-attainment
was previously established in
\cite{HamamotoTakahashi2021JFA}.

\begin{cor}[Exact quadratic transition and remainder identities]
\label{cor:quadratic-transition}
Let \(p=2\), \(\gamma\ne0\), and let \(\psi\in\mathcal A_{2,\gamma}\)
have \(\mathcal J_{2,\gamma}[\psi]<\infty\). Recall that
\[
t=\log r,\qquad b=\gamma-1,
\qquad V(t,\theta)=e^{bt}\psi(e^t e_r(\theta)),
\]
and put
\[
c=b^2+4b+1=(\gamma+1)^2-3.
\]
Separate the angular mean by writing
\[
V=\overline V+\varphi,
\qquad
\overline V(t)=\frac1{2\pi}\int_{\Sone}V(t,\theta)\,\dd\theta,
\qquad \int_{\Sone}\varphi(t,\theta)\,\dd\theta=0.
\]
Set
\[
f=\overline V_t-b\overline V,\qquad
P=f+\varphi_t-b\varphi,\qquad Q=\varphi_\theta,
\qquad U=(P,Q).
\]
All unqualified integrals below are over \(\R\times\Sone\), including
those involving \(f\), which depends only on \(t\). Define
\begin{align*}
W[\varphi]&=\int(\varphi_\theta^2-\varphi^2),\\
G_b[\varphi]&=\int
\left[(\varphi_{t\theta}^2-\varphi_t^2)
+\varphi_{\theta\theta}^2+(c-1)\varphi_\theta^2-c\varphi^2\right].
\end{align*}
Both are nonnegative and vanish on the first angular mode
\(\varphi=A(t)\cos\theta+B(t)\sin\theta\).

Write \(\mathcal I=\mathcal I_{2,\gamma}\),
\(\mathcal J=\mathcal J_{2,\gamma}\), and
\[
H_\gamma=1+\gamma^2+\frac{c}{b^2+1}
=\frac{\gamma^2(b^2+3)}{b^2+1}.
\]
The exact deficit identities are
\begin{align}
\mathcal J-(1+\gamma^2)\mathcal I
&=\int(|U_t|^2+\varphi_t^2+c\varphi^2)
+G_b[\varphi],&&c\ge0,
\label{eq:quadratic-outer}\displaybreak[0]\\
\mathcal J-H_\gamma\mathcal I
&=\int|U_t|^2
+\frac{-c\int f^2-4b\int\varphi_t^2}{b^2+1}
+G_b[\varphi]-\frac{c}{b^2+1}W[\varphi],&&c<0.
\label{eq:quadratic-inner}
\end{align}
In each stated region, every term on the right is nonnegative.
Consequently,
\[
C(2,\gamma)=
\begin{cases}
\displaystyle\frac1{1+\gamma^2},&c\ge0,\\[2mm]
\displaystyle\frac{b^2+1}{\gamma^2(b^2+3)},&c<0,
\end{cases}
\]
in agreement with \eqref{eq:Mazya-p2}. The transition values are
\(\gamma=-1\pm\sqrt3\).
\end{cor}

\begin{rem}[Comparison with the established quadratic remainder]
In the notation of \cite[Theorem~1]{HamamotoTakahashi2021JFA}, their
radial derivative remainder is
\(\int|\partial_t Z|^2\,\dd t\,\dd\theta\) in our logarithmic
cylinder. Since \(Z=P e_r+Q e_\theta\) and the polar basis is
independent of \(t\), this equals \(\int|U_t|^2\), the common term in
\eqref{eq:quadratic-outer} and \eqref{eq:quadratic-inner}.
The remaining angular contributions are organized here through
\(W[\varphi]\) and \(G_b[\varphi]\). We do not claim term-by-term
domination of the angular remainder in
\cite{HamamotoTakahashi2021JFA}.
\end{rem}

\begin{rem}[Non-attainment and modal rigidity]
\label{rem:quadratic-rigidity}
Both deficit identities contain \(\int|U_t|^2\).
Equality would force \(U\) to be independent of \(t\).
Since \(U\in L^2(\R\times\Sone)\), this implies \(U=0\),
and hence \(\nabla\psi=0\). Compact support then gives \(\psi=0\).
Thus no nonzero admissible function attains either sharp branch.

For \(c>0\), the outer deficit controls
\[
\|U_t\|_2^2+\|\varphi_t\|_2^2+c\|\varphi\|_2^2.
\]
For \(c<0\), let \(\Pi_1\) be the angular projection onto
\(\operatorname{span}\{\cos\theta,\sin\theta\}\).
The Fourier representation of \(W\) gives
\(W[\varphi]\ge3\|\varphi-\Pi_1\varphi\|_2^2\), so the inner
deficit controls
\[
\|U_t\|_2^2
+\frac{-c}{b^2+1}\|f\|_2^2
+\frac{-4b}{b^2+1}\|\varphi_t\|_2^2
+\frac{-3c}{b^2+1}\|\varphi-\Pi_1\varphi\|_2^2.
\]
Thus, for normalized almost extremizers, the outer region suppresses
the angular part, while the inner region suppresses the radial
contribution and higher angular modes. In both cases, the surviving
fields vary slowly in \(t\).
At \(c=0\), the radial and first angular modes have the same sharp
quotient, so the deficit no longer selects between them. Higher
angular modes remain penalized by \(G_b\).
These are rigidity statements for almost extremizing sequences,
without a claim of convergence to a nonzero extremizer.
\end{rem}

\subsection{The critical logarithmic endpoint}

The divergence of $C(p,\gamma)$ in
Theorem~\ref{thm:asymptotic} leads us to a logarithmic
replacement at $\gamma=1-2/p$. For $p>1$, the classical
Hardy ground-state transform gives the sharp inequality
and an exact deficit identity. The argument does not
use the curl-free constraint and follows the convexity
method of
\cite{FrankSeiringer2008,CazacuFlynnLamLu2024,DoFlynnLamLu2023}.

\begin{thm}[Sharp logarithmic endpoint and exact deficit]
\label{thm:critical-log-identity}
Let \(1<p<\infty\), \(R>0\), and
\(\psi\in C_0^\infty(B_R;\R)\). Then
\begin{equation}\label{eq:critical-log}
\int_{B_R}\frac{|\nabla\psi(x)|^p}
{|x|^2\bigl(\log(R/|x|)\bigr)^p}\,\dd x
\le\left(\frac{p}{p-1}\right)^p
\int_{B_R}|D^2\psi(x)|^p|x|^{p-2}\,\dd x.
\end{equation}
The constant is sharp even for radial functions and is not attained
by any nonzero admissible \(\psi\).

To describe the deficit, write
\begin{align*}
\mathcal H[\psi]&:=\int_{B_R}
\frac{|\nabla\psi|^p}{|x|^2\bigl(\log(R/|x|)\bigr)^p}\,\dd x,\\
\mathcal E[\psi]&:=\int_{B_R}|D^2\psi|^p|x|^{p-2}\,\dd x,
\qquad
\mathcal D[\psi]:=\mathcal E[\psi]-\alpha^p\mathcal H[\psi],
\end{align*}
where \(\alpha=(p-1)/p\). Set
\[
t=\log(R/r)>0,\qquad s=\log t,
\qquad
w(s,\theta)=t^{-\alpha}\nabla\psi(Re^{-t}e_r(\theta)),
\quad t=e^s.
\]
The field \(w\) is expressed in fixed Cartesian components.
For
\[
X=(\alpha w+w_s,e^s w_\theta),
\qquad
Y=(w_s,e^s w_\theta),
\]
viewed as vectors in \(\R^4\), the deficit has the exact representation
\begin{equation}\label{eq:critical-log-identity}
\mathcal D[\psi]
=\int_{\R\times\Sone}\Cp(X,Y)\,\dd s\,\dd\theta.
\end{equation}
Consequently, for constants \(c_p>0\) depending only on \(p\),
\begin{align}
\mathcal D[\psi]
&\ge c_p\int_{\R\times\Sone}
\bigl(|w_s|^2+e^{2s}|w_\theta|^2\bigr)^{p/2},
&&p\ge2,\label{eq:critical-log-stability-large}\\
\mathcal D[\psi]
&\ge c_p\frac{\left[\displaystyle\int_{\R\times\Sone}
\bigl(|w_s|^2+e^{2s}|w_\theta|^2\bigr)^{p/2}\right]^{2/p}}
{\bigl(\mathcal E[\psi]+\alpha^p\mathcal H[\psi]\bigr)^{(2-p)/p}},
&&1<p<2.\label{eq:critical-log-stability-small}
\end{align}
For \(\psi=0\), the right-hand side of
\eqref{eq:critical-log-stability-small} is understood as zero.
In particular, for \(p=2\),
\begin{equation}\label{eq:critical-log-p2-identity}
\mathcal E[\psi]-\frac14\mathcal H[\psi]
=\int_{\R\times\Sone}
\bigl(|w_s|^2+e^{2s}|w_\theta|^2\bigr)\,\dd s\,\dd\theta.
\end{equation}
\end{thm}

\begin{rem}[Meaning of stability at the endpoint]
If \(\mathcal H[\psi_n]=1\) and \(\mathcal D[\psi_n]\to0\), then
\[
\|w_{n,s}\|_p+\|e^s w_{n,\theta}\|_p\to0.
\]
The estimates give the rate \(O(\mathcal D^{1/p})\) for \(p\ge2\)
and \(O(\mathcal D^{1/2})\) for \(1<p<2\).
On a bounded \(s\)-interval, \(e^s\) is bounded below, so
Poincar\'e's inequality controls the distance from the average vector
on that cylinder by these derivative norms. The constant depends on
the interval. There is no nonzero global \(L^p\) extremizer:
after normalization, the sharpness sequence constructed in the proof of Theorem~\ref{thm:critical-log-identity} moves towards
\(s=-\infty\) while varying increasingly slowly.
The conclusion is therefore control of transformed derivatives,
without a claim of convergence to a fixed nonzero extremizer.
\end{rem}

\begin{rem}[The \(p=2\) endpoint]
For \(p=2\), the sharp inequality is
\(\frac14\mathcal H\le\mathcal E\), equivalently
\(\mathcal H\le4\mathcal E\). This agrees with
Hamamoto and Takahashi \cite{HamamotoTakahashi2021Critical}.
\end{rem}

\begin{rem}[Relation with the one-dimensional Hardy identity]
For comparison, let \(v\in C_0^\infty((0,\infty))\) be scalar and
put \(g(t)=t^\alpha v(t)\). The same ground-state calculation gives
\[
\int_0^\infty
\left(|g'|^p-\alpha^p\frac{|g|^p}{t^p}\right)\,\dd t
=\int_0^\infty\frac1t
\Cp(\alpha v+tv',tv')\,\dd t.
\]
Indeed, \(g'=t^{\alpha-1}(\alpha v+tv')\), and the linear term in
the remainder expansion integrates to
\(\alpha^{p-1}[|v|^p]_0^\infty=0\).
For a radial potential, the full Hessian energy also contains the
angular contribution:
\[
|D^2\psi|^2=r^{-2}\bigl(|g'|^2+|g|^2\bigr).
\]
The identity \eqref{eq:critical-log-identity} retains this contribution
and applies to all admissible scalar potentials.
Its vector-valued proof uses the Euclidean convexity estimates of
Lemma~\ref{lem:Cp}, without the curl-free constraint.
The Bregman expansion is the same general mechanism
as in nonlinear ground-state representations
\cite{FrankSeiringer2008}.
General weighted $L^p$ identities based on related
convexity remainders were developed by Do, Flynn,
Lam and Lu \cite{DoFlynnLamLu2023}, while a closely
related quadratic Caffarelli--Kohn--Nirenberg
identity framework was established in
\cite{CazacuFlynnLamLu2024}.
See also
\cite{ShaimerdenovYessirkegenovZhangirbayev2026}.
\end{rem}

The assumption \(p>1\) is essential, as the following radial
counterexample shows.

\begin{prop}[Failure of the critical logarithmic estimate at \(p=1\)]
\label{prop:p1-critical-failure}
For every \(R>0\), there is no finite constant \(C\) such that
\begin{equation}\label{eq:p1-critical-fail}
\int_{B_R}
\frac{|\nabla\psi(x)|}{|x|^2\log\frac{R}{|x|}}\,\dd x
\le C\int_{B_R}|D^2\psi(x)|\,|x|^{-1}\,\dd x
\end{equation}
holds for every \(\psi\in C_0^\infty(B_R)\).
\end{prop}

\section{Proofs}\label{sec:all-proofs}
\label{SEC 4}

We prove the auxiliary statements first and then establish the main results in the order described below.

\subsection{Preliminary lemmas and propositions}

\begin{proof}[Proof of Lemma \ref{lem:Cp}]
The nonnegativity is the Bregman-divergence inequality for \(F(z)=|z|^p\).

Assume first \(p\ge2\). Put \(x=\xi-\eta\) and \(y=\eta\). If \(x\neq0\), write
\[
|x+y|^2=|x|^2(1+s),
\qquad
s=\frac{2x\cdot y+|y|^2}{|x|^2}\ge-1.
\]
Since \((1+s)^{p/2}\ge1+(p/2)s\),
\[
|x+y|^p
\ge
|x|^p+p|x|^{p-2}x\cdot y
+\frac p2|x|^{p-2}|y|^2,
\]
which is \eqref{eq:Cp-quadratic-explicit}. The case \(x=0\) is immediate. The
estimate \eqref{eq:Cp-coercive-large} follows from the monotonicity inequality
\[
\bigl(|z|^{p-2}z-|x|^{p-2}x\bigr)\cdot(z-x)
\ge2^{2-p}|z-x|^p,
\]
see \cite[Chapter~12, inequality~(I)]{Lindqvist2019}. Taking \(z=x+ty\),
dividing by \(t\), multiplying by \(p\), and integrating over
\(0<t<1\) gives
\[
|x+y|^p-|x|^p-p|x|^{p-2}x\cdot y
\ge p2^{2-p}|y|^p\int_0^1t^{p-1}\,\dd t
=2^{2-p}|y|^p.
\]

Now let \(1<p<2\). For \(F(z)=|z|^p\), the smallest eigenvalue of \(D^2F(z)\) is
\(p(p-1)|z|^{p-2}\). Taylor's formula with integral remainder, first applied to a
regularization of \(F\) if the segment crosses the origin, gives
\[
\begin{aligned}
\Cp(\xi,\eta)
&\ge
p(p-1)|\eta|^2
\int_0^1(1-t)|\xi-\eta+t\eta|^{p-2}\,\dd t
\\
&\ge
\frac{p(p-1)}2
\bigl(|\xi-\eta|+|\eta|\bigr)^{p-2}|\eta|^2,
\end{aligned}
\]
which is \eqref{eq:Cp-subquadratic-explicit}. Since
\[
|\xi-\eta|+|\eta|
\le 2\bigl(|\xi|+|\xi-\eta|\bigr),
\]
this gives \eqref{eq:Cp-coercive-small} with the displayed constant.
\end{proof}

\begin{proof}[Proof of Lemma \ref{lem:1d-amplitude}]
For \(1<p<\infty\), apply Lemma~\ref{lem:Cp} with
\(\xi=c'-\sigma c\), \(\eta=c'\). Since \(\xi-\eta=-\sigma c\),
\[
\Cp(c'-\sigma c,c')
=
|c'-\sigma c|^p-|c|^p
+p\sigma |c|^{p-2}c\cdot c'.
\]
Now \(|c|^p\in W^{1,1}(\R)\), hence it has zero limits at
\(\pm\infty\), and
\[
\int_\R p|c|^{p-2}c\cdot c'\,\dd s
=
\int_\R (|c|^p)'\,\dd s=0.
\]
This proves \eqref{eq:1d-amplitude-remainder}.

For \(p=1\), the Euclidean norm is convex. At almost every point where \(c\neq0\),
\[
|c'-\sigma c|
\ge
-\sigma\frac{c}{|c|}\cdot(c'-\sigma c)
=
|c|-\sigma(|c|)'.
\]
The chain rule for \(|c|\in W^{1,1}(\R)\) and integration over \(\R\) give
\eqref{eq:1d-amplitude}. On \(\{c=0\}\), one has \(c'=0\) almost everywhere, so the same
inequality holds there with the zero subgradient.
\end{proof}

\begin{proof}[Proof of Lemma \ref{lem:cylinder}]
Since \(r=e^t\), the area element becomes
\[
\dd x=r\,\dd r\,\dd\theta
=r^2\,\dd t\,\dd\theta.
\]
Using \(\nabla\psi=r^{-a}Z\), we obtain
\[
\begin{aligned}
|\nabla\psi|^p r^{p(\gamma-1)}\,\dd x
&=
|Z|^p r^{-pa+p(\gamma-1)+2}
\,\dd t\,\dd\theta\\
&=
|Z|^p\,\dd t\,\dd\theta,
\end{aligned}
\]
because \(pa=p(\gamma-1)+2\).
Integration proves \eqref{eq:cylinder-I}.

Next, let us write
\(\psi_{ij}=\partial_{x_i}\partial_{x_j}\psi\).
All derivatives below are evaluated at
\(x=(r\cos\theta,r\sin\theta)\).
The Hessian and its Frobenius norm are
\[
D^2\psi=
\begin{pmatrix}
\psi_{11}&\psi_{12}\\
\psi_{12}&\psi_{22}
\end{pmatrix},
\qquad
|D^2\psi|^2
=\psi_{11}^2+2\psi_{12}^2+\psi_{22}^2.
\]

Differentiating \(Z=r^a\nabla\psi\) by the chain rule gives
\[
\begin{aligned}
Z_t-aZ
&=
r^{a+1}
\begin{pmatrix}
\psi_{11}\cos\theta+\psi_{12}\sin\theta\\
\psi_{12}\cos\theta+\psi_{22}\sin\theta
\end{pmatrix}
=
r^{a+1}D^2\psi\,e_r,
\\[1ex]
Z_\theta
&=
r^{a+1}
\begin{pmatrix}
-\psi_{11}\sin\theta+\psi_{12}\cos\theta\\
-\psi_{12}\sin\theta+\psi_{22}\cos\theta
\end{pmatrix}
=
r^{a+1}D^2\psi\,e_\theta.
\end{aligned}
\]
Here \(aZ\) is the term produced by differentiating \(r^a\).

Squaring the four components and adding cancels the mixed
products. Thus
\[
\begin{aligned}
|Z_t-aZ|^2+|Z_\theta|^2
&=
r^{2(a+1)}
\left(\psi_{11}^2+2\psi_{12}^2+\psi_{22}^2\right)\\
&=
r^{2(a+1)}|D^2\psi|^2.
\end{aligned}
\]
Consequently,
\[
\begin{aligned}
|D^2\psi|^p r^{p\gamma}\,\dd x
&=
\left(|Z_t-aZ|^2+|Z_\theta|^2\right)^{p/2}\\
&\qquad{}\times
r^{-p(a+1)+p\gamma+2}\,\dd t\,\dd\theta.
\end{aligned}
\]
The remaining power of \(r\) is zero, since
\(p(a+1)=p\gamma+2\).
Integration proves \eqref{eq:cylinder-J}.

The Hessian is symmetric, so
\[
e_r\cdot D^2\psi\,e_\theta
=
e_\theta\cdot D^2\psi\,e_r.
\]
Using the two chain-rule formulas above,
\[
\begin{aligned}
e_r\cdot Z_\theta
&=
r^{a+1}e_r\cdot D^2\psi\,e_\theta\\
&=
r^{a+1}e_\theta\cdot D^2\psi\,e_r\\
&=
e_\theta\cdot(Z_t-aZ).
\end{aligned}
\]
This proves \eqref{eq:cylinder-constraint}.
\end{proof}

\begin{proof}[Proof of Proposition \ref{prop:first-harmonic}]
Take \(\phi(\theta)=\cos\theta\). Since
\(\phi'=-\sin\theta\) and \(\phi''-b\phi=-a\cos\theta\),
\[
\begin{aligned}
a^2b^2\phi^2+2a^2(\phi')^2+(\phi''-b\phi)^2
&=a^2\bigl(b^2\cos^2\theta+1+\sin^2\theta\bigr),\\
b^2\phi^2+(\phi')^2&=b^2\cos^2\theta+\sin^2\theta.
\end{aligned}
\]
Thus \(\mathcal Q_{p,a}[\cos\theta]=|a|^pJ_p(b)/I_p(b)\), and
\eqref{eq:angular-test-bound} gives the assertions. The quantities
\(I_p,J_p\) are the angular integrals for this test, not the weighted
energies \(\Ipg,\Jpg\).

For \(p=2\), integration of \(\sin^2\theta\) and \(\cos^2\theta\)
gives \(I_2(b)=\pi(b^2+1)\) and \(J_2(b)=\pi(b^2+3)\). Hence
\[
\mathcal Q_{2,a}[\cos\theta]-(1+a^2)
=\frac{2a^2-b^2-1}{b^2+1}
=\frac{b^2+4b+1}{b^2+1}.
\]
Since \(a=\gamma\) and \(b=\gamma-1\), condition
\eqref{eq:first-harmonic-breaking} is equivalent to
\[
-1-\sqrt3<\gamma<-1+\sqrt3,\qquad \gamma\ne0.
\]
\end{proof}

\begin{proof}[Proof of Lemma \ref{lem:KL-identity}]
Using integration by parts and \(P_\theta=Q_t-bQ\), we obtain
\[
\begin{aligned}
K
&=-2\int |U|^{p-2}QP_\theta
  -\int PQ\,\partial_\theta(|U|^{p-2})\\
&=2bA+\int Q^2\partial_t(|U|^{p-2})
  -\int PQ\,\partial_\theta(|U|^{p-2})\\
&=2bA+(p-2)\int |U|^{p-4}
  \left[Q^2(U\cdot U_t)-PQ(U\cdot U_\theta)\right].
\end{aligned}
\]
The second line uses
\(\int\partial_t(|U|^{p-2}Q^2)=0\).
Expanding the determinants and using \(Q_t-P_\theta=bQ\) gives
\[
\begin{aligned}
Q\left[Q\det(U,U_\theta)+P\det(U,U_t)\right]
&=PQ(U\cdot U_\theta)-Q^2(U\cdot U_t)\\
&\quad+b|U|^2Q^2.
\end{aligned}
\]
Thus the last integral in the first calculation is \(bA-L\), and
\[
K=2bA+(p-2)(bA-L)=pbA-(p-2)L.
\]

For completeness, finite energy gives \(U,U_t,U_\theta\in L^p\).
The integrand of \(L\) is bounded in absolute value by
\[
|U|^{p-1}\bigl(|U_t|^2+|U_\theta|^2\bigr)^{1/2},
\]
so all the displayed integrals are finite by H\"older's inequality.
The differential calculations are first written where \(U\ne0\).
At zeros, they are justified by applying the chain rule to
\(|U|^{p-2}PQ\) and \(|U|^{p-2}Q^2\), which extend as \(C^1\)
functions of \(U\) for \(p>1\). Periodicity removes the angular
boundary terms. Compact cutoffs in \(t\), followed by a limit, remove
its boundary terms because these two functions belong to \(W^{1,1}\).
\end{proof}

\begin{proof}[Proof of Lemma \ref{lem:exact-deficit}]
By \eqref{eq:U-cylinder-J},
\(\Jpg[\psi]=\int|X+Y|^p\), while
\(|X|^p=\lambda^p|U|^p\). The definition of \(\Cp\) gives
\[
|X+Y|^p-|X|^p
=p|X|^{p-2}X\cdot Y+\Cp(X+Y,Y).
\]
Since
\[
\int |X|^{p-2}X\cdot Y
=\lambda^{p-2}\left[
-a\int |U|^{p-2}U\cdot U_t+K\right]
=\lambda^{p-2}K,
\]
where \(\int |U|^{p-2}U\cdot U_t=p^{-1}\int\partial_t|U|^p=0\),
integration proves the identity. Nonnegativity of \(R_p\) follows
from Lemma~\ref{lem:Cp}.
\end{proof}

\begin{proof}[Proof of Lemma \ref{lem:angular-cell}]
First, the proposed field satisfies the two constraints. Indeed,
\[
\mathcal Y[q]
=
\begin{pmatrix}
-\sin 2\theta & \cos 2\theta\\
 \cos 2\theta & \sin 2\theta
\end{pmatrix}q,
\]
so its Cartesian components integrate to zero. Moreover,
orthogonality of \(e_r,e_\theta\) gives
\[
e_r\cdot\mathcal Y[q]=e_\theta\cdot q,
\qquad
|\mathcal Y[q]|^2
=(e_\theta\cdot q)^2+(e_r\cdot q)^2
=|q|^2.
\]
Consequently, its energy is
\[
\int_{\Sone}
\left(|q|^2+|\mathcal Y[q]|^2\right)^{p/2}\,\dd\theta
=2\pi\,2^{p/2}|q|^p.
\]

Next, \(Y\) and \(\mathcal Y[q]\) have the same radial component,
so their difference is tangential:
\[
Y-\mathcal Y[q]=h(\theta)e_\theta
\]
for some scalar function \(h\). Their zero averages imply
\[
0=\int_{\Sone}h e_\theta\,\dd\theta
=
\left(
-\int_{\Sone}h\sin\theta\,\dd\theta,\,
 \int_{\Sone}h\cos\theta\,\dd\theta
\right).
\]
Therefore,
\[
\int_{\Sone}h e_r\,\dd\theta=0,
\]
and hence
\[
\begin{aligned}
\int_{\Sone}\mathcal Y[q]\cdot
       (Y-\mathcal Y[q])\,\dd\theta
&=\int_{\Sone}h(e_r\cdot q)\,\dd\theta\\
&=q\cdot\int_{\Sone}h e_r\,\dd\theta=0.
\end{aligned}
\]

If \(q=0\), both conclusions follow directly. Assume \(q\ne0\).
For \(1<p<\infty\), the definition of \(\Cp\) gives
\[
\begin{aligned}
\left(|q|^2+|Y|^2\right)^{p/2}
&=2^{p/2}|q|^p\\
&\quad+p(\sqrt2\,|q|)^{p-2}
 \mathcal Y[q]\cdot(Y-\mathcal Y[q])\\
&\quad+\Cp\bigl((q,Y),(0,Y-\mathcal Y[q])\bigr).
\end{aligned}
\]
Integration removes the middle term and proves
\eqref{eq:angular-cell-remainder}. Since \(\Cp\ge0\), it also
proves \eqref{eq:angular-cell-lower}.

For \(p=1\), Cauchy--Schwarz gives
\[
\begin{aligned}
\sqrt{|q|^2+|Y|^2}
&\ge
\frac{|q|^2+\mathcal Y[q]\cdot Y}{\sqrt2\,|q|}\\
&=\sqrt2\,|q|
+\frac{\mathcal Y[q]\cdot(Y-\mathcal Y[q])}
       {\sqrt2\,|q|}.
\end{aligned}
\]
The last term again integrates to zero, proving the required bound.
\end{proof}

\begin{proof}[Proof of Proposition \ref{prop:constrained-liminf}]
We show that every sequence \(\varepsilon_n\downarrow0\) and
\(W_n\in\mathscr X_{\varepsilon_n,\sigma}^p\) with
\(\|W_n\|_{L^p}=1\) and uniformly bounded energy satisfies the
required lower bound. This suffices by choosing fields whose
energies exceed the infima by at most \(1/n\) along a sequence
realizing the lower limit of \(m_{p,\sigma}\). If that lower
limit is infinite, the conclusion is immediate.

Below, unqualified integrals and \(L^p\)-norms are over
\(\R\times\Sone\).

\medskip
\noindent
\emph{Step 1. Smooth the fields without increasing the energy.}

Put
\[
d_n=(W_n)_s-\sigma W_n,
\qquad
Y_n=\varepsilon_n^{-1}(W_n)_\theta.
\]
Bounded energy controls \(d_n,Y_n\) in \(L^p\), and
\[
(W_n)_s=d_n+\sigma W_n
\]
is bounded there as well.

Choose \(0\le\rho\in C_0^\infty((-1,1))\) with
\(\int_\R\rho=1\), and set
\[
\delta_n=\sqrt{\varepsilon_n},
\qquad
\rho_{\delta_n}(s)=\delta_n^{-1}\rho(s/\delta_n).
\]
Average only in \(s\):
\[
\widetilde W_n=\rho_{\delta_n}*_sW_n,
\qquad
\widetilde d_n=\rho_{\delta_n}*_sd_n,
\qquad
\widetilde Y_n=\rho_{\delta_n}*_sY_n.
\]
Averaging commutes with derivatives. Since \(e_r,e_\theta\)
are independent of \(s\), we retain
\[
\widetilde d_n
=(\widetilde W_n)_s-\sigma\widetilde W_n,
\qquad
\widetilde Y_n
=\varepsilon_n^{-1}(\widetilde W_n)_\theta,
\]
and the constraint is preserved exactly:
\[
e_r\cdot\widetilde Y_n
=e_\theta\cdot\widetilde d_n.
\]
Jensen's inequality gives
\begin{equation}\label{eq:mollification-energy}
\int
\left(|\widetilde d_n|^2+|\widetilde Y_n|^2\right)^{p/2}
\le
\int\left(|d_n|^2+|Y_n|^2\right)^{p/2}.
\end{equation}
The translation estimate also gives
\[
\|\widetilde W_n-W_n\|_{L^p}
\le C\delta_n\|(W_n)_s\|_{L^p}
=O(\sqrt{\varepsilon_n}),
\]
and consequently
\[
\|\widetilde W_n\|_{L^p}^p=1+o(1).
\]
These estimates hold also for \(p=1\).

\medskip
\noindent
\emph{Step 2. Make the radial expression nearly independent
of the angle.}

Since \((W_n)_\theta=\varepsilon_nY_n\),
\[
(d_n)_\theta
=\varepsilon_n\bigl((Y_n)_s-\sigma Y_n\bigr)
\]
in distributions. After averaging, the derivative in \(s\)
falls on the kernel:
\[
(\widetilde Y_n)_s
=(\rho_{\delta_n}')*_sY_n,
\qquad
\|\rho_{\delta_n}'\|_{L^1(\R)}
=\delta_n^{-1}\|\rho'\|_{L^1(\R)}.
\]
Thus no bound for \((Y_n)_s\) is needed, and
\[
\|(\widetilde d_n)_\theta\|_{L^p}
\le C\varepsilon_n(\delta_n^{-1}+1)
=O(\sqrt{\varepsilon_n}).
\]

Define the angular average
\[
d_n^0(s)
=\frac1{2\pi}\int_{\Sone}
\widetilde d_n(s,\theta)\,\dd\theta.
\]
Periodic Poincar\'e gives
\begin{equation}\label{eq:dn-d0}
\|\widetilde d_n-d_n^0\|_{L^p}
=O(\sqrt{\varepsilon_n})=o(1).
\end{equation}
Moreover, Jensen's inequality gives
\[
2\pi\int_\R|d_n^0|^p\,\dd s
\le\int|\widetilde d_n|^p
\le\int|d_n|^p\le C.
\]
Both statements remain valid at \(p=1\).

\medskip
\noindent
\emph{Step 3. Apply the angular calculation with a small error.}

For each fixed \(s\), set
\[
Y_n^0=\mathcal Y[d_n^0],
\qquad
H_n=\widetilde Y_n-Y_n^0.
\]
Periodicity gives
\(\int_{\Sone}\widetilde Y_n\,\dd\theta=0\), while
Lemma~\ref{lem:angular-cell} gives the same for \(Y_n^0\).
Hence
\[
\int_{\Sone}H_n\,\dd\theta=0,
\qquad
e_r\cdot H_n
=e_\theta\cdot(\widetilde d_n-d_n^0).
\]

Let \(J(v_1,v_2)=(-v_2,v_1)\). The formula for the angular
cell gives
\[
Y_n^0
=Jd_n^0+2(e_\theta\cdot d_n^0)e_r.
\]
Since \(d_n^0\) is independent of \(\theta\), the term
\(Jd_n^0\) has zero integral against \(H_n\). Therefore,
\begin{equation}\label{eq:linear-cell-error}
\begin{aligned}
\int_{\Sone}Y_n^0\cdot H_n\,\dd\theta
&=
2\int_{\Sone}(e_\theta\cdot d_n^0)
\bigl(e_\theta\cdot(\widetilde d_n-d_n^0)\bigr)
\,\dd\theta,\\
\left|
\int_{\Sone}Y_n^0\cdot H_n\,\dd\theta
\right|
&\le
2|d_n^0|
\int_{\Sone}|\widetilde d_n-d_n^0|\,\dd\theta.
\end{aligned}
\end{equation}

For \(p>1\), convexity at \((d_n^0,Y_n^0)\), whose length
is \(\sqrt2|d_n^0|\), gives, where \(d_n^0\ne0\),
\[
\begin{aligned}
\left(|\widetilde d_n|^2+
|\widetilde Y_n|^2\right)^{p/2}
&\ge 2^{p/2}|d_n^0|^p\\
&\quad+
p(\sqrt2|d_n^0|)^{p-2}
\left[
d_n^0\cdot(\widetilde d_n-d_n^0)
+Y_n^0\cdot H_n
\right].
\end{aligned}
\]
The first linear term has zero angular integral.
Using \eqref{eq:linear-cell-error} for the second gives
\[
\begin{aligned}
&\int_{\Sone}
\left(|\widetilde d_n|^2+
|\widetilde Y_n|^2\right)^{p/2}\,\dd\theta\\
&\quad\ge
2\pi\,2^{p/2}|d_n^0|^p
-C_p|d_n^0|^{p-1}
\int_{\Sone}|\widetilde d_n-d_n^0|\,\dd\theta.
\end{aligned}
\]
Where \(d_n^0=0\), this bound follows from nonnegativity.
After integration in \(s\), H\"older's inequality bounds
the error by
\[
C_p
\left(2\pi\int_\R|d_n^0|^p\,\dd s\right)^{(p-1)/p}
\|\widetilde d_n-d_n^0\|_{L^p}
=o(1).
\]

For \(p=1\) and \(d_n^0\ne0\), Cauchy--Schwarz gives
\[
\begin{aligned}
\sqrt{|\widetilde d_n|^2+|\widetilde Y_n|^2}
&\ge \sqrt2|d_n^0|\\
&\quad+
\frac{
d_n^0\cdot(\widetilde d_n-d_n^0)+Y_n^0\cdot H_n
}{\sqrt2|d_n^0|}.
\end{aligned}
\]
Angular integration therefore yields
\[
\begin{aligned}
&\int_{\Sone}
\sqrt{|\widetilde d_n|^2+|\widetilde Y_n|^2}
\,\dd\theta\\
&\quad\ge
2\pi\sqrt2|d_n^0|
-\sqrt2\int_{\Sone}
|\widetilde d_n-d_n^0|\,\dd\theta.
\end{aligned}
\]
If \(d_n^0=0\), the right-hand side is nonpositive,
so the same bound is immediate. Its integrated error
tends to zero by \eqref{eq:dn-d0}.

Thus, for every \(1\le p<\infty\),
\begin{equation}\label{eq:liminf-cell-bound}
\int
\left(|\widetilde d_n|^2+|\widetilde Y_n|^2\right)^{p/2}
\ge
2\pi\,2^{p/2}\int_\R|d_n^0|^p\,\dd s-o(1).
\end{equation}

\medskip
\noindent
\emph{Step 4. Use the one-dimensional inequality
and the normalization.}

Let
\[
c_n(s)=\frac1{2\pi}\int_{\Sone}
\widetilde W_n(s,\theta)\,\dd\theta.
\]
Then \(c_n\in W^{1,p}(\R;\R^2)\) and
\[
d_n^0=c_n'-\sigma c_n.
\]
Lemma~\ref{lem:1d-amplitude} gives
\[
\int_\R|d_n^0|^p\,\dd s
\ge\int_\R|c_n|^p\,\dd s.
\]

Since
\((\widetilde W_n)_\theta
=\varepsilon_n\widetilde Y_n\),
periodic Poincar\'e also gives
\[
\|\widetilde W_n-c_n\|_{L^p}
\le C\varepsilon_n\|\widetilde Y_n\|_{L^p}
=O(\varepsilon_n).
\]
Together with
\(\|\widetilde W_n-W_n\|_{L^p}=o(1)\) and
\(\|W_n\|_{L^p}=1\), this implies
\[
2\pi\int_\R|c_n|^p\,\dd s=1+o(1).
\]

Combining the estimates,
\[
\begin{aligned}
\int\left(|d_n|^2+|Y_n|^2\right)^{p/2}
&\ge
\int\left(
|\widetilde d_n|^2+|\widetilde Y_n|^2
\right)^{p/2}\\
&\ge
2\pi\,2^{p/2}\int_\R|d_n^0|^p\,\dd s-o(1)\\
&\ge
2\pi\,2^{p/2}\int_\R|c_n|^p\,\dd s-o(1)\\
&=2^{p/2}-o(1).
\end{aligned}
\]
Applying this to the almost-minimizing fields described
at the beginning proves \eqref{eq:constrained-liminf}.
\end{proof}

\begin{proof}[Proof of Lemma \ref{lem:L-control}]
All integrals below are over \(\R\times\Sone\). Recall that
\(U=(P,Q)\) satisfies
\[
Q_t-P_\theta=bQ.
\]
Put \(r=|U|\), and recall
\[
A=\int r^{p-2}Q^2,
\qquad
L=\int r^{p-4}Q
\left[Q\det(U,U_\theta)+P\det(U,U_t)\right].
\]
We estimate \(L\) and \(L-bA\) using, respectively, the changes
in the direction and length of \(U\).

First work on \(\{r>0\}\). Since
\[
rr_t=PP_t+QQ_t,
\qquad
rr_\theta=PP_\theta+QQ_\theta,
\]
the elementary planar identity gives
\[
|U_t|^2
=r_t^2+\frac{\det(U,U_t)^2}{r^2},
\qquad
|U_\theta|^2
=r_\theta^2+\frac{\det(U,U_\theta)^2}{r^2}.
\]
Consequently,
\[
\begin{aligned}
B={}&\int r^{p-2}(r_t^2+r_\theta^2)\\
&+\int r^{p-4}
\left[\det(U,U_t)^2+\det(U,U_\theta)^2\right].
\end{aligned}
\]
The first integral measures changes in length, while the second
measures changes in direction.

By Cauchy--Schwarz,
\[
\left|Q\det(U,U_\theta)+P\det(U,U_t)\right|
\le r
\left[\det(U,U_t)^2+\det(U,U_\theta)^2\right]^{1/2}.
\]
Applying Cauchy--Schwarz once more, now to the integral
defining \(L\), yields
\[
L^2\le
A\int r^{p-4}
\left[\det(U,U_t)^2+\det(U,U_\theta)^2\right].
\]

To estimate \(L-bA\), expand the determinants and use the
constraint:
\[
\begin{aligned}
&Q\det(U,U_\theta)+P\det(U,U_t)
-r(Pr_\theta-Qr_t)\\
&\qquad=r^2(Q_t-P_\theta)=br^2Q.
\end{aligned}
\]
Multiplying by \(r^{p-4}Q\) and integrating gives
\[
L-bA=\int r^{p-3}Q(Pr_\theta-Qr_t).
\]
Since
\[
|Pr_\theta-Qr_t|
\le r(r_t^2+r_\theta^2)^{1/2},
\]
Cauchy--Schwarz gives
\[
(L-bA)^2
\le
\left(\int r^{p-2}Q^2\right)
\left(\int r^{p-2}(r_t^2+r_\theta^2)\right)
=
A\int r^{p-2}(r_t^2+r_\theta^2).
\]
Adding the two estimates and using the decomposition of \(B\)
proves \eqref{eq:strong-L-control}.

The integrals are finite by H\"older's inequality, since
\(U,U_t,U_\theta\in L^p\) and \(p\ge2\).
Expressions involving negative powers of \(r\) are assigned
value zero on \(\{U=0\}\). This is consistent because
\(U_t\) and \(U_\theta\) vanish almost everywhere on that set.
\end{proof}

\subsection{Near-critical asymptotics and rigidity}

\begin{proof}[Proof of Theorem \ref{thm:asymptotic}]
Let
\[
\Lambda_p(a):=\frac1{C(p,\gamma)},
\qquad
a=\gamma-1+\frac2p.
\]
Set \(\varepsilon=|a|\), \(\sigma=\operatorname{sgn}(a)\), and
\(s=\varepsilon t\). If \(Z\) is the field from Lemma~\ref{lem:cylinder}, define
\[
W(s,\theta)=\varepsilon^{-1/p}Z(s/\varepsilon,\theta).
\]
Then
\[
\int|W|^p\,\dd s\,\dd\theta
=
\int|Z|^p\,\dd t\,\dd\theta,
\]
and
\[
\frac1{\varepsilon^p}
\int
\left(|Z_t-aZ|^2+|Z_\theta|^2\right)^{p/2}
\,\dd t\,\dd\theta
=
\int
\left(
|W_s-\sigma W|^2+
\left|\frac{W_\theta}{\varepsilon}\right|^2
\right)^{p/2}
\,\dd s\,\dd\theta.
\]
Moreover, \(W\in\mathscr X_{\varepsilon,\sigma}^p\). Hence
\[
\frac{\Lambda_p(a)}{|a|^p}
\ge
m_{p,\sigma}(\varepsilon),
\]
and Proposition~\ref{prop:constrained-liminf} gives
\begin{equation}\label{eq:lambda-liminf}
\liminf_{a\to0}\frac{\Lambda_p(a)}{|a|^p}
\ge2^{p/2}.
\end{equation}

The matching upper bound follows from Proposition~\ref{prop:first-harmonic}:
\begin{equation}\label{eq:first-mode-upper}
\Lambda_p(a)\le |a|^p\frac{J_p(b)}{I_p(b)}.
\end{equation}
As \(a\to0\), \(b\to-1\), and
\[
I_p(-1)=2\pi,
\qquad
J_p(-1)=2\pi\,2^{p/2}.
\]
Thus
\[
\limsup_{a\to0}\frac{\Lambda_p(a)}{|a|^p}
\le2^{p/2}.
\]
Together with \eqref{eq:lambda-liminf}, this proves
\eqref{eq:main-asymptotic}.
\end{proof}

\begin{proof}[Proof of Theorem \ref{thm:rigidity}]
The estimates in Proposition~\ref{prop:constrained-liminf} give
\[
\|\widetilde W_n-W_n\|_{L^p}
=O(\sqrt{\varepsilon_n}),
\qquad
\|\widetilde d_n-d_n^0\|_{L^p}
=O(\sqrt{\varepsilon_n}),
\]
and
\[
\|\widetilde W_n-c_n\|_{L^p}=O(\varepsilon_n).
\]
Since \(\|W_n\|_{L^p}=1\), these prove
\eqref{eq:rig-mass}, \eqref{eq:rig-W}, and \eqref{eq:rig-d}.

To obtain the remaining conclusions, retain the two nonnegative
convexity remainders in the lower-bound proof. Set
\[
X_n=(\widetilde d_n,\widetilde Y_n),
\qquad
X_n^0=(d_n^0,\mathcal Y[d_n^0]).
\]
The angular expansion used in
Proposition~\ref{prop:constrained-liminf} gives
\begin{equation}\label{eq:rig-angular-decomp}
\int_{\R\times\Sone}|X_n|^p
=
2\pi\,2^{p/2}\int_\R|d_n^0|^p\,\dd s
+R_n^{\rm ang}+o(1),
\end{equation}
where
\[
R_n^{\rm ang}
:=
\int_{\R\times\Sone}
\Cp(X_n,X_n-X_n^0)\,\dd s\,\dd\theta
\ge0.
\]
The \(o(1)\) term follows from the linear-error estimate
\eqref{eq:linear-cell-error}, H\"older's inequality, and
\(\widetilde d_n-d_n^0\to0\) in \(L^p\).

Likewise, Lemma~\ref{lem:1d-amplitude} gives
\begin{equation}\label{eq:rig-1d-decomp}
\int_\R|d_n^0|^p\,\dd s
=
\int_\R|c_n|^p\,\dd s+R_n^{\rm 1D},
\qquad
R_n^{\rm 1D}:=
\int_\R\Cp(d_n^0,c_n')\,\dd s\ge0.
\end{equation}

Mollification does not increase the energy. Applying
Proposition~\ref{prop:constrained-liminf} to the normalized
mollified fields, separately for the two possible signs
\(\sigma_n\), gives the matching lower bound. Hence
\[
\int_{\R\times\Sone}|X_n|^p=2^{p/2}+o(1).
\]
Combining \eqref{eq:rig-angular-decomp} and
\eqref{eq:rig-1d-decomp}, we obtain
\[
\begin{aligned}
\int_{\R\times\Sone}|X_n|^p
={}&2\pi\,2^{p/2}\int_\R|c_n|^p\,\dd s\\
&+2\pi\,2^{p/2}R_n^{\rm 1D}
+R_n^{\rm ang}+o(1).
\end{aligned}
\]
The left-hand side and the first term on the right both tend
to \(2^{p/2}\). Nonnegativity therefore implies
\[
R_n^{\rm ang}\longrightarrow0,
\qquad
R_n^{\rm 1D}\longrightarrow0.
\]

For \(p\ge2\), estimate \eqref{eq:Cp-coercive-large} gives
\[
\|X_n-X_n^0\|_{L^p}\longrightarrow0,
\qquad
\|c_n'\|_{L^p(\R)}\longrightarrow0.
\]
For \(1<p<2\), the same conclusions follow from
\eqref{eq:Cp-coercive-small}, the uniform \(L^p\)-bounds,
and Remark~\ref{rem:p-less-2}.
This proves \eqref{eq:rig-flat} and \eqref{eq:rig-Y}.

Finally,
\[
d_n^0+\sigma_nc_n=c_n',
\qquad
\mathcal Y[d_n^0]+\sigma_n\mathcal Y[c_n]
=\mathcal Y[c_n'].
\]
The last two assertions follow from these identities,
\eqref{eq:rig-d}, \eqref{eq:rig-Y}, and
\(\|c_n'\|_{L^p(\R)}\to0\).
\end{proof}

\subsection{Radial sharpness and symmetry breaking}

\begin{proof}[Proof of Theorem \ref{thm:radial}]
Let \(\psi(x)=f(r)\), \(r=|x|>0\). 
Since
$\nabla\psi=f'(r)\frac{x}{|x|},$ we have $|\nabla\psi|=|f'|$.
Differentiation of $\nabla\psi$ gives
\begin{align*}
D^2\psi
&=
\begin{pmatrix}
\dfrac{\partial^2\psi}{\partial x_1^2}
&
\dfrac{\partial^2\psi}{\partial x_1\partial x_2}
\\[1ex]
\dfrac{\partial^2\psi}{\partial x_2\partial x_1}
&
\dfrac{\partial^2\psi}{\partial x_2^2}
\end{pmatrix}
=
\begin{pmatrix}
\dfrac{x_1^2}{r^2}f''(r)+\dfrac{x_2^2}{r^3}f'(r)
&
\dfrac{x_1x_2}{r^2}f''(r)-\dfrac{x_1x_2}{r^3}f'(r)
\\[2ex]
\dfrac{x_1x_2}{r^2}f''(r)-\dfrac{x_1x_2}{r^3}f'(r)
&
\dfrac{x_2^2}{r^2}f''(r)+\dfrac{x_1^2}{r^3}f'(r)
\end{pmatrix}.
\end{align*}
Consequently, using \(x_1^2+x_2^2=r^2\),
\[
\begin{aligned}
|D^2\psi|^2
&=
\left(\frac{x_1^2}{r^2}f''
      +\frac{x_2^2}{r^3}f'\right)^2
+2\left(\frac{x_1x_2}{r^2}f''
       -\frac{x_1x_2}{r^3}f'\right)^2+
\left(\frac{x_2^2}{r^2}f''
      +\frac{x_1^2}{r^3}f'\right)^2\\
&=
\frac{(f'')^2}{r^4}(x_1^2+x_2^2)^2
+\frac{(f')^2}{r^6}(x_1^2+x_2^2)^2+
\frac{2f''f'}{r^5}
\left(x_1^2x_2^2-2x_1^2x_2^2+x_1^2x_2^2\right)\\
&=(f'')^2+\frac{(f')^2}{r^2}.
\end{aligned}
\]
Set
\[
a=\gamma-1+\frac2p,\qquad r=e^t,
\qquad v(t)=e^{at}f'(e^t).
\]
Thus
\[
f'(e^t)=e^{-at}v,
\qquad
f''(e^t)=e^{-(a+1)t}(v'-av),
\]
and
\[
\begin{aligned}
\frac{\Ipg[\psi]}{2\pi}
&=
\int_\R |f'(e^t)|^p
e^{[p(\gamma-1)+2]t}\,\dd t\\
&=\int_\R |v|^p\,\dd t,
\\[1ex]
\frac{\Jpg[\psi]}{2\pi}
&=
\int_\R
\left(|f''(e^t)|^2+e^{-2t}|f'(e^t)|^2\right)^{p/2}
e^{(p\gamma+2)t}\,\dd t\\
&=
\int_\R\left((v'-av)^2+v^2\right)^{p/2}\,\dd t.
\end{aligned}
\]
Finite energy gives \(v,v'-av\in L^p(\R)\), hence
\(v\in W^{1,p}(\R)\). Consequently,
\[
|v|^p\in W^{1,1}(\R),
\qquad
\int_\R (|v|^p)'\,\dd t=0.
\]

For \(1<p<\infty\), the definition of \(\Cp\) gives
\[
\begin{aligned}
\Cp\bigl((v'-av,v),(v',0)\bigr)
&=
\left((v'-av)^2+v^2\right)^{p/2}
-(1+a^2)^{p/2}|v|^p\\
&\quad+pa(1+a^2)^{(p-2)/2}|v|^{p-2}vv'\\
&=
\left((v'-av)^2+v^2\right)^{p/2}
-(1+a^2)^{p/2}|v|^p\\
&\quad+a(1+a^2)^{(p-2)/2}(|v|^p)'.
\end{aligned}
\]
Integrating, we obtain the exact identity
\[
\Jpg[\psi]-(1+a^2)^{p/2}\Ipg[\psi]
=
2\pi\int_\R
\Cp\bigl((v'-av,v),(v',0)\bigr)\,\dd t
\ge0.
\]

For \(p=1\), Cauchy--Schwarz gives, on \(\{v\neq0\}\),
\[
\begin{aligned}
\sqrt{(v'-av)^2+v^2}
&\ge
\frac{-av(v'-av)+v^2}{\sqrt{1+a^2}\,|v|}\\
&=
\sqrt{1+a^2}\,|v|
-\frac{a}{\sqrt{1+a^2}}(|v|)'.
\end{aligned}
\]
On \(\{v=0\}\), the same inequality holds almost everywhere
because \((|v|)'=0\) there. Integrating and using
\(\int_\R(|v|)'\,\dd t=0\), we obtain
\[
\Jpg[\psi]\ge\sqrt{1+a^2}\,\Ipg[\psi].
\]

For sharpness, fix \(0\neq v\in C_0^\infty(\R)\) and define
\[
f_L(e^t)
:=-\int_t^\infty e^{(1-a)s}v(s/L)\,\dd s,
\qquad
\psi_L(x):=f_L(|x|).
\]
Then
\[
e^{at}f_L'(e^t)=v(t/L).
\]
The function \(\psi_L\) is smooth, constant near the origin,
and vanishes outside a ball. Hence
\[
\psi_L\in C_0^\infty(\R^2),
\qquad
\nabla\psi_L(0)=0,
\]
and both weighted integrals are finite.

Rescaling the integration variable gives
\[
\begin{aligned}
\frac{\Jpg[\psi_L]}{\Ipg[\psi_L]}
&=
\frac{
\displaystyle\int_\R
\left((L^{-1}v'-av)^2+v^2\right)^{p/2}\,\dd t
}{
\displaystyle\int_\R |v|^p\,\dd t
}\\
&\longrightarrow (1+a^2)^{p/2}
\qquad (L\to\infty).
\end{aligned}
\]
The limit follows by dominated convergence, since the
integrands are uniformly bounded and supported in a
fixed compact set.

Taking the infimum over admissible nonzero radial functions
with finite weighted integrals, the preceding inequality
and the test sequence yield
\[
\begin{aligned}
(1+a^2)^{p/2}
&\le
\inf_{\psi}
\frac{\Jpg[\psi]}{\Ipg[\psi]}\\
&\le
\frac{\Jpg[\psi_L]}{\Ipg[\psi_L]}
\longrightarrow (1+a^2)^{p/2}.
\end{aligned}
\]
Therefore,
\[
C_{\rm rad}(p,\gamma)^{-1}
=
\inf_{\psi}\frac{\Jpg[\psi]}{\Ipg[\psi]}
=
(1+a^2)^{p/2},
\]
and hence
\[
C_{\rm rad}(p,\gamma)=(1+a^2)^{-p/2}.
\]
Although the full inequality fails at the critical parameter
\(\gamma=1-2/p\), its restriction to radial functions remains
valid, with sharp constant \(C_{\rm rad}=1\).
\end{proof}

\begin{proof}[Proof of Theorem \ref{thm:radial-sharpness}]
Put \(\lambda=\sqrt{1+a^2}\). By Lemma~\ref{lem:exact-deficit},
\eqref{eq:outer-Cp-deficit} holds with \(X=(-aU,JU)\) and
\(Y=(U_t,U_\theta)\).

Assume first \(p\ge2\). By \eqref{eq:Cp-quadratic-explicit},
\[
R_p
\ge
\frac p2\lambda^{p-2}B.
\]
Thus Lemmas~\ref{lem:KL-identity} and \ref{lem:L-control} give
\[
\begin{aligned}
\Jpg-\lambda^p\Ipg
&\ge
p\lambda^{p-2}\left[pbA-(p-2)L+\frac12B\right]
\\
&\ge
p\lambda^{p-2}
\left[
pbA-(p-2)L+
\frac{L^2+(L-bA)^2}{2A}
\right]
\end{aligned}
\]
when \(A>0\). If \(A=0\), then \(L=K=0\) and the conclusion is immediate.
For \(A>0\), put \(x=L/A\). The bracket divided by \(A\) is
\[
f_{p,b}(x)
=
pb-(p-2)x+\frac12\bigl(x^2+(x-b)^2\bigr).
\]
Its minimum is
\[
\min_{x\in\R}f_{p,b}(x)
=
\frac14\left[b^2+2(p+2)b-(p-2)^2\right].
\]
Therefore \(\Jpg\ge\lambda^p\Ipg\) whenever
\[
b^2+2(p+2)b-(p-2)^2\ge0.
\]
Its two branches are exactly \eqref{eq:radial-region-pge2-b}.

Now let \(1<p<2\) and put \(q=2-p\). By
\eqref{eq:Cp-subquadratic-explicit},
\[
R_p
\ge
\kappa_p\int
\frac{|Y|^2}{(|X|+|Y|)^q}.
\]
Write
\[
u:=|U|,
\qquad
v:=\left(|U_t|^2+|U_\theta|^2\right)^{1/2}.
\]
Since \(|X|=\lambda u\),
\[
R_p
\ge
\kappa_p\int\frac{v^2}{(\lambda u+v)^q}.
\]
Moreover, Lemma~\ref{lem:KL-identity} gives
\[
K=pbA+qL,
\]
and, from the definition of \(L\),
\[
L\ge -\int u^{p-2}|Q|v.
\]
On \(\{u>0\}\), define
\[
s:=\frac{|Q|}{u}\in[0,1],
\qquad
r:=\frac{v}{\lambda u}\ge0.
\]
Then \eqref{eq:outer-Cp-deficit} yields
\[
\Jpg-\lambda^p\Ipg
\ge
\lambda^{p-2}\int_{\{u>0\}}u^p
\left[
p^2bs^2-pq\lambda rs
+\kappa_p\lambda^2\frac{r^2}{(1+r)^q}
\right].
\]
On \(\{u=0\}\), the \(A\)- and \(K\)-terms vanish and the \(\Cp\)-remainder is
nonnegative. If \(s=0\), the displayed bracket is nonnegative. If \(s>0\), write
\(r=sz\). Since \(1+sz\le1+z\),
\[
\frac{1}{(1+sz)^q}\ge\frac{1}{(1+z)^q},
\]
so the bracket divided by \(s^2\) is bounded below by
\[
p^2b-pq\lambda z
+\kappa_p\lambda^2\frac{z^2}{(1+z)^q}.
\]
Condition \eqref{eq:radial-region-sub2} therefore implies
\(\Jpg\ge\lambda^p\Ipg\).

In either case,
\[
C(p,\gamma)\le\lambda^{-p}=(1+a^2)^{-p/2}.
\]
The reverse inequality follows from Theorem~\ref{thm:radial}, because radial functions
belong to the full admissible class. Hence equality holds.
\end{proof}

\begin{proof}[Proof of Corollary \ref{cor:radial-sharpness-sub2}]
Let \(q=2-p\). For \(0\le z\le1\),
\[
\frac{z^2}{(1+z)^q}\ge2^{p-2}z^2,
\]
whereas for \(z\ge1\),
\[
\frac{z^2}{(1+z)^q}\ge2^{p-2}z^p.
\]
Using these bounds in \eqref{eq:radial-region-sub2} and maximizing
the resulting quadratic on $[0,1]$ and $Az-Dz^p$ on
$[1,\infty)$ gives the two bounds below. In the second one,
the remaining factor $\lambda^{(p-2)/(p-1)}$ is at most one,
since $\lambda\ge1$ and $p<2$. Thus it is sufficient to have
\[
b\ge
\max\left\{
\frac{(2-p)^2}{4d_p},
\frac{p-1}{p^2}(2-p)^{\frac{p}{p-1}}d_p^{-\frac1{p-1}}
\right\}
=B_p.
\]
Since \(b=\gamma-2+2/p\), this is \eqref{eq:radial-region-sub2-explicit}.
\end{proof}

\begin{proof}[Proof of Theorem \ref{thm:negative-sub2}]
Put \(q=2-p\), \(\kappa=p(p-1)/2\), \(d=-b\), and
\(\lambda=\sqrt{1+a^2}=\sqrt{1+(1-d)^2}\). Our hypothesis gives
\(d\ge D_p>0\). Write
\(u=|U|\) and \(v=(|U_t|^2+|U_\theta|^2)^{1/2}\).
The subquadratic estimate \eqref{eq:Cp-subquadratic-explicit} and
Lemma~\ref{lem:KL-identity} imply
\[
R_p\ge\kappa\int\frac{v^2}{(\lambda u+v)^q},
\qquad K=-pdA+qL,
\qquad L\ge-\int u^{p-2}|Q|v.
\]
The last bound follows by Cauchy--Schwarz on the two determinants
in the definition of \(L\). On \(\{u>0\}\), set
\(s=|Q|/u\in[0,1]\) and \(z=v/u\ge0\).
The deficit identity therefore gives
\begin{equation}\label{eq:negative-bracket}
\Jpg-\lambda^p\Ipg
\ge\lambda^{-q}\int_{\{u>0\}}u^p
\left[\frac{\kappa z^2}{(1+z/\lambda)^q}
-p^2ds^2-pqsz\right].
\end{equation}
The omitted contribution on \(\{u=0\}\) is nonnegative. The
curl-free constraint gives the additional pointwise bound
\[
ds=\frac{|P_\theta-Q_t|}{u}
\le\frac{\sqrt2v}{u}=\sqrt2z,
\qquad s\le\min\{1,\sqrt2z/d\}.
\]
For \(0\le z\le d/\sqrt2\), the bracket in
\eqref{eq:negative-bracket} is bounded below by
\[
z^2\left[\frac{\kappa}{(1+z/\lambda)^q}
-\frac{2p^2+\sqrt2pq}{d}\right].
\]
Since \(2\lambda^2-d^2=(d-2)^2\ge0\), we have
\(z/\lambda\le1\). The displayed expression is nonnegative when
\(d\ge D_p=2^q(2p^2+\sqrt2pq)/\kappa\).

For \(z\ge d/\sqrt2\), it suffices to show
\[
\frac{\kappa z^2}{(1+z/\lambda)^q}\ge p^2d+pqz.
\]
The quotient of the left side by the right side is increasing in
\(z>0\), because its logarithmic derivative exceeds
\((1-q)/z>0\). At \(z=d/\sqrt2\) this inequality is the first case.
Thus \(\Jpg\ge\lambda^p\Ipg\) for every admissible potential.
The radial test sequence of Theorem~\ref{thm:radial} gives
the reverse bound on the optimal constant.
\end{proof}

\begin{proof}[Proof of Proposition \ref{prop:instability}]
Take \(\phi_\varepsilon=1+\varepsilon\cos\theta\) and
\(\lambda=\sqrt{1+a^2}\). Since \(S_p(1)=p/2>0\), the
hypothesis implies \(b\ne0\). Both densities can therefore be
expanded uniformly for sufficiently small \(\varepsilon\). Write $c=\cos\theta$ and $s=\sin\theta$. The denominator
and numerator densities for $\phi_\varepsilon$ are respectively
\[
b^2+2b^2\varepsilon c+\varepsilon^2(b^2c^2+s^2)
\quad\text{and}\quad
\lambda^2b^2+2b(a^2b+a)\varepsilon c
+a^2\varepsilon^2\bigl((b^2+1)c^2+2s^2\bigr).
\]
For $b\ne0$, expand their $p/2$ powers around the constant terms
and use $\int c=0$, $\int c^2=\int s^2=\pi$. The quotient is
\[
\mathcal Q_{p,a}[\phi_\varepsilon]
=\lambda^p+\frac{p\lambda^{p-2}}{4b^2}
S_p(a)\varepsilon^2+o(\varepsilon^2).
\]
For \(S_p(a)<0\), choosing a sufficiently small nonzero
\(\varepsilon\) proves \eqref{eq:strict-breaking}.
\end{proof}

\subsection{The quadratic case}

\begin{proof}[Proof of Corollary \ref{cor:quadratic-transition}]
\emph{The basic quadratic identity.}
At \(p=2\), the Bregman remainder in \eqref{eq:outer-Cp-deficit}
is \(\int(|U_t|^2+|U_\theta|^2)\), while
\eqref{eq:KL-identity} gives \(K=2b\int Q^2\). Hence
\begin{equation}\label{eq:quadratic-basic}
\mathcal J-(1+\gamma^2)\mathcal I
=\int(|U_t|^2+|U_\theta|^2+4bQ^2).
\end{equation}

The cylinder representation gives \(U,U_t,U_\theta\in L^2\).
Since \(Q=\varphi_\theta\) and
\(Q_t=\varphi_{t\theta}\), the zero angular mean and the circle
Poincar\'e inequality imply \(\varphi,\varphi_t\in L^2\).
Also, \(\varphi_{\theta\theta}=Q_\theta\in L^2\).
These facts justify the integrations by parts and Fourier expansions
below, without requiring \(\overline V\in L^2\).

\emph{Nonnegative angular remainders.}
The circle Poincar\'e inequality immediately gives
\(W[\varphi]\ge0\). To examine \(G_b\), write
\[
\varphi(t,\theta)=\sum_{n\ne0}\varphi_n(t)e^{in\theta}.
\]
Parseval's identity gives
\[
G_b[\varphi]
=2\pi\sum_{n\ne0}(n^2-1)
\int_\R
\left(|\varphi_n'|^2+(n^2+c)|\varphi_n|^2\right)\,\dd t.
\]
Since \(c=(b+2)^2-3\ge-3\), every summand is nonnegative:
the factor \(n^2-1\) is zero for \(|n|=1\), and
\(n^2+c\ge1\) for \(|n|\ge2\).
This proves \(G_b[\varphi]\ge0\) and explains why both remainders
vanish on the first angular mode.

\emph{The outer deficit.}
From \(P_\theta=\varphi_{t\theta}-b\varphi_\theta\) and
\(Q_\theta=\varphi_{\theta\theta}\), we obtain
$$
\int(|U_\theta|^2+4bQ^2)
=\int\left[
\varphi_{t\theta}^2+\varphi_{\theta\theta}^2
+(b^2+4b)\varphi_\theta^2\right]
=G_b[\varphi]+\int(\varphi_t^2+c\varphi^2).
$$
The cross term vanishes because
\(\int\varphi_{t\theta}\varphi_\theta
=\frac12\int\partial_t(\varphi_\theta^2)=0\).
Substitution into \eqref{eq:quadratic-basic} proves
\eqref{eq:quadratic-outer}. This algebraic identity holds for every
\(c\). When \(c\ge0\), all its terms are nonnegative, giving
\(\mathcal J\ge(1+\gamma^2)\mathcal I\).

\emph{The inner deficit.}
Angular orthogonality gives
\[
\mathcal I=\int(P^2+Q^2)
=\int\left[f^2+(\varphi_t-b\varphi)^2
+\varphi_\theta^2\right].
\]
Since \(\int\varphi_t\varphi=0\), this becomes
\[
\mathcal I
=\int(f^2+\varphi_t^2+b^2\varphi^2+\varphi_\theta^2)
=\int(f^2+\varphi_t^2)+(b^2+1)\int\varphi^2+W[\varphi].
\]
Subtracting \(c\mathcal I/(b^2+1)\) from the outer identity cancels
\(c\int\varphi^2\). The coefficient of \(\int\varphi_t^2\) becomes
\[
1-\frac{c}{b^2+1}=-\frac{4b}{b^2+1},
\]
which proves \eqref{eq:quadratic-inner}.
If \(c<0\), then \(-2-\sqrt3<b<-2+\sqrt3<0\), so every term
on its right is nonnegative. Thus
\(\mathcal J\ge H_\gamma\mathcal I\).

\emph{Sharpness and the transition.}
For \(c\ge0\), the radial cutoff sequence from
Theorem~\ref{thm:radial} approaches the quotient \(1+\gamma^2\),
proving sharpness of the outer branch.
For \(c<0\), take
\[
V_L(t,\theta)=\eta(t/L)\cos\theta,
\qquad 0\ne\eta\in C_0^\infty(\R),\qquad L\to\infty.
\]
Here \(\overline V_L=0\), so \(f=0\) and \(\varphi=V_L\).
Since this is a first angular mode,
\(W[\varphi]=G_b[\varphi]=0\). Direct integration gives
\[
\mathcal I
=\pi(b^2+1)L\int_\R\eta^2\,\dd t
+\frac\pi L\int_\R|\eta'|^2\,\dd t
\asymp L.
\]
Moreover, \(\int\varphi_t^2=O(L^{-1})\) and
\(\int|U_t|^2=O(L^{-1})\). The inner identity therefore yields
\[
\mathcal J-H_\gamma\mathcal I=O(L^{-1}),
\qquad
\frac{\mathcal J}{\mathcal I}
=H_\gamma+O(L^{-2})\longrightarrow H_\gamma.
\]
The corresponding potentials
\(\psi_L(re_r(\theta))=r^{-b}\eta((\log r)/L)\cos\theta\)
are smooth and supported in annuli, so they are admissible.
This proves sharpness of the inner branch.
Finally, \(H_\gamma-(1+\gamma^2)=c/(b^2+1)\), so the branches
cross exactly at \(c=0\), or \(\gamma=-1\pm\sqrt3\).
\end{proof}

\subsection{The logarithmic endpoint}

\begin{proof}[Proof of Theorem \ref{thm:critical-log-identity}]
\emph{Change of variables.}
Put
$u(t,\theta)=\nabla\psi(Re^{-t}e_r(\theta))
$
in fixed Cartesian components. Since \(\dd t/\dd r=-1/r\),
\[
\partial_r\nabla\psi=-r^{-1}u_t,
\qquad
r^{-1}\partial_\theta\nabla\psi=r^{-1}u_\theta.
\]
These are the radial and tangential columns of the Hessian, so its
Frobenius norm satisfies
\[
|D^2\psi|^2=r^{-2}\bigl(|u_t|^2+|u_\theta|^2\bigr).
\]
Using \(\dd x=r\,\dd r\,\dd\theta\) and
\(\dd r/r=-\dd t\), we obtain
\[
\mathcal H
=\int_0^\infty\int_{\Sone}t^{-p}|u|^p\,\dd\theta\,\dd t,
\qquad
\mathcal E
=\int_0^\infty\int_{\Sone}
\bigl(|u_t|^2+|u_\theta|^2\bigr)^{p/2}\,\dd\theta\,\dd t.
\]
Now \(u(t,\theta)=t^\alpha w(\log t,\theta)\), and hence
$
u_t=t^{\alpha-1}(\alpha w+w_s)$ and $u_\theta=t^\alpha w_\theta$. The identities \(\alpha p=p-1\) and \(\dd t/t=\dd s\) therefore give
\[
\mathcal H=\int_{\R\times\Sone}|w|^p\,\dd s\,\dd\theta,
\qquad
\mathcal E=\int_{\R\times\Sone}
\bigl(|\alpha w+w_s|^2+e^{2s}|w_\theta|^2\bigr)^{p/2}
\,\dd s\,\dd\theta.
\]

Both quantities are finite directly from the assumptions on \(\psi\).
Indeed, \(u\) vanishes for sufficiently small \(t\) and is bounded
as \(t\to\infty\), while \(\int_1^\infty t^{-p}\,\dd t<\infty\).
Also, \(D^2\psi\) is bounded and the radial factor in \(\mathcal E\)
is \(r^{p-1}\), which is integrable at zero.
Moreover,
\[
w_s+\alpha w=t^{1-\alpha}u_t,
\qquad
\int_{\R\times\Sone}|w_s+\alpha w|^p\,\dd s\,\dd\theta
=\int_0^\infty\int_{\Sone}|u_t|^p\,\dd\theta\,\dd t
\le\mathcal E.
\]
Thus \(w,w_s\in L^p\), which justifies the integration by parts below.

\emph{Exact deficit and the inequality.}
Recall that
$
X=(\alpha w+w_s,e^s w_\theta)$ and 
$
Y=(w_s,e^s w_\theta)$ are vectors in \(\R^4\). Thus \(X-Y=(\alpha w,0)\), and the
definition of \(\Cp\) yields
\[
\begin{aligned}
\Cp(X,Y)
&=|X|^p-\alpha^p|w|^p
-p\alpha^{p-1}|w|^{p-2}w\cdot w_s\\
&=|X|^p-\alpha^p|w|^p
-\alpha^{p-1}\partial_s|w|^p.
\end{aligned}
\]
Here \(|w|^p\) and \(\partial_s|w|^p\) are integrable by
H\"older's inequality. Integrating the derivative with cutoffs in \(s\)
therefore gives zero. This proves \eqref{eq:critical-log-identity}.
Since \(\Cp\ge0\), it follows that
\(\mathcal E\ge\alpha^p\mathcal H\), which is
\eqref{eq:critical-log}.

\emph{Derivative estimates and nonattainment.}
For \(p\ge2\), Lemma~\ref{lem:Cp} gives
\(\Cp(X,Y)\ge c_p|Y|^p\), proving
\eqref{eq:critical-log-stability-large}.
For \(1<p<2\), the same lemma gives
\[
\mathcal D\ge c_p\int
\frac{|Y|^2}{(|X|+|X-Y|)^{2-p}}.
\]
H\"older's inequality, with exponents \(2/p\) and \(2/(2-p)\), implies
\[
\|Y\|_p^p
\le C_p\mathcal D^{p/2}
\left[\int(|X|+|X-Y|)^p\right]^{(2-p)/2}
\le C_p\mathcal D^{p/2}
(\mathcal E+\alpha^p\mathcal H)^{(2-p)/2}.
\]
Rearrangement proves \eqref{eq:critical-log-stability-small}.
At \(p=2\), \(\mathscr C_2(X,Y)=|Y|^2\), giving
\eqref{eq:critical-log-p2-identity}.

If equality holds in \eqref{eq:critical-log}, these estimates force
\(Y=0\), so \(w_s=0\) and \(w_\theta=0\).
The field \(w\) is therefore constant. Since
\(w\in L^p(\R\times\Sone)\), this constant is zero.
Hence \(\nabla\psi=0\), and compact support gives \(\psi=0\).

\emph{Sharpness in the radial class.}
Fix \(0\ne f\in C_0^\infty((-1,1))\), and define
\[
g_n(t)=t^\alpha f\left(\frac{\log t}{n}+n\right),
\qquad
\psi_n(r)=-\int_r^R
g_n\left(\log\frac{R}{\rho}\right)\,\dd\rho.
\]
For each \(n\), \(g_n\in C_0^\infty((0,\infty))\).
Thus \(\psi_n\) is zero near \(r=R\) and constant near \(r=0\),
so it defines a function in \(C_0^\infty(B_R)\).
Its gradient and transformed field are
\[
\nabla\psi_n(re_r)=g_n(\log(R/r))e_r,
\qquad
w_n(s,\theta)=f\left(\frac{s}{n}+n\right)e_r(\theta).
\]
Since \(\partial_\theta e_r=e_\theta\), it follows that
\[
(w_n)_s=\frac1n f'\left(\frac{s}{n}+n\right)e_r,
\qquad
(w_n)_\theta=f\left(\frac{s}{n}+n\right)e_\theta.
\]
With \(y=s/n+n\), the factors \(2\pi n\) cancel in the quotient,
and we obtain
\[
\frac{\mathcal E[\psi_n]}{\mathcal H[\psi_n]}
=
\frac{\displaystyle\int_{-1}^{1}
\left[
\left|\alpha f(y)+\frac{f'(y)}{n}\right|^2
+e^{2n(y-n)}|f(y)|^2
\right]^{p/2}\,\dd y}
{\displaystyle\int_{-1}^{1}|f(y)|^p\,\dd y}.
\]
Here \(f'/n\to0\), while
\(e^{2n(y-n)}\le e^{-2n^2+2n}\to0\) on \([-1,1]\).
Dominated convergence therefore gives
\[
\frac{\mathcal E[\psi_n]}{\mathcal H[\psi_n]}
\longrightarrow\alpha^p.
\]
The stretching makes the \(s\)-derivative small, and the translation
towards \(s=-\infty\) makes the weighted angular derivative small.
Thus \(\alpha^{-p}=(p/(p-1))^p\) is the optimal constant in
\eqref{eq:critical-log}.
\end{proof}

\begin{proof}[Proof of Proposition \ref{prop:p1-critical-failure}]
For \(0<\varepsilon<1/4\), choose
\(g_\varepsilon\in C_0^\infty((0,\infty))\) satisfying
\[
0\le g_\varepsilon\le1,
\qquad
g_\varepsilon=1\text{ on }[2\varepsilon,1],
\qquad
\operatorname{supp}g_\varepsilon\subset[\varepsilon,2].
\]
Choose the rise from zero to one and the fall back to zero to be
monotone. Then
\[
\int_0^\infty|g_\varepsilon'|\,\dd t=2,
\qquad
\int_0^\infty g_\varepsilon\,\dd t\le2.
\]
Define the radial potential
\[
\psi_\varepsilon(r)
=-\int_r^R g_\varepsilon\left(\log\frac{R}{\rho}\right)\,\dd\rho.
\]
It is constant near the origin and zero near the boundary, and hence
belongs to \(C_0^\infty(B_R)\).
The radial change of variables gives
\[
\int_{B_R}
\frac{|\nabla\psi_\varepsilon|}{|x|^2\log(R/|x|)}\,\dd x
=2\pi\int_0^\infty\frac{g_\varepsilon(t)}t\,\dd t
\ge2\pi\log\frac1{2\varepsilon}
\longrightarrow\infty.
\]
On the other hand,
\begin{align*}
\int_{B_R}|D^2\psi_\varepsilon|\,|x|^{-1}\,\dd x
&=2\pi\int_0^\infty
\sqrt{|g_\varepsilon'|^2+|g_\varepsilon|^2}\,\dd t\\
&\le2\pi\left(\int_0^\infty|g_\varepsilon'|\,\dd t
+\int_0^\infty g_\varepsilon\,\dd t\right)
\le8\pi.
\end{align*}
Thus the left-hand side of \eqref{eq:p1-critical-fail} is unbounded
while the right-hand energy remains bounded, proving the claim.
\end{proof}

\section*{Acknowledgements}
AI tools were used in the preparation of this manuscript. All mathematical arguments and proofs presented in the paper were developed and verified by the authors.

\end{document}